\documentclass[10pt, reqno]{amsart}

\usepackage{hyperref} 
\usepackage{amssymb, mathrsfs,bbold}
\usepackage{amsmath, amsthm}
\usepackage{enumerate}
\usepackage{dsfont}
\usepackage{color}
\usepackage[a4paper]{geometry}
\usepackage[utf8]{inputenc}   
\usepackage{stmaryrd}
\usepackage{titlesec, wasysym}
\usepackage{appendix}
\usepackage{todonotes, fancyhdr}
\usepackage{chngcntr}
\usepackage{cancel}

\usepackage{tikz-cd}
\usetikzlibrary{calc}

\usepackage{setspace}
\renewcommand{\baselinestretch}{0.99}

\usepackage[all]{xy}
\numberwithin{subsection}{section}
\numberwithin{subsubsection}{subsection}
\numberwithin{equation}{section} 

\newenvironment{Dem}[1][\unskip]{%
    \begin{list}{\hspace{1.15cm}{\sf \textbf{{\small Proof #1 --}}}}{%
        \setlength{\topsep}{0pt}%
        \setlength{\leftmargin}{0pt}%
        \setlength{\rightmargin}{0pt}%
        \setlength{\listparindent}{0pt}%
        \setlength{\itemindent}{0pt}%
        \setlength{\parsep}{0pt}%
        \addtolength{\leftmargin}{20pt}%
        \addtolength{\rightmargin}{0pt}%
    } \item }{\hfill $\rhd$\end{list}\smallskip}

\newenvironment{Dem*}[1][\unskip]{%
    \begin{list}{\hspace{0cm}{\sf \textbf{{\small Proof #1 --}}}}{%
        \setlength{\topsep}{0pt}%
        \setlength{\leftmargin}{0pt}%
        \setlength{\rightmargin}{0pt}%
        \setlength{\listparindent}{0pt}%
        \setlength{\itemindent}{0pt}%
        \setlength{\parsep}{0pt}%
        \addtolength{\leftmargin}{20pt}%
        \addtolength{\rightmargin}{0pt}%
    } \item }{\hfill $\rhd$\end{list}\smallskip}

\renewcommand\thesection       {\arabic{section}}
\renewcommand\thesubsection    {\thesection{\boldmath $.$}\arabic{subsection}}
\renewcommand\thesubsubsection    {\thesection{\boldmath $.$}\arabic{subsection}{\boldmath $.$}\arabic{subsubsection}}

\titleformat{\section}[block]
{\filcenter\normalfont\sffamily\bfseries\Large}
{{\hspace{-0.7cm}}\thesection \hspace{0.2em} --\vspace{0.3cm}}{0.5em}{}

\titleformat{\subsection}[block]
{\filcenter\normalfont\sffamily\bfseries\large}  						  
{\hspace{-0.7cm}\thesubsection \hspace{0.5em} \vspace{0.3cm}}{.5em}{}  
\titlespacing{\subsection}{-0pc}{1.5ex plus .1ex minus .2ex}{0pc}

\titleformat{\subsubsection}[block]
{\normalfont\sffamily\bfseries}					  
{\thesubsubsection \vspace{0.3cm}}{.5em}{}  
\titlespacing{\subsection}{-0pc}{1.5ex plus .1ex minus .2ex}{0pc}

\newtheoremstyle{mystyle}
{3pt}               
{3pt}               
{\it }                      
{}                      
{\sffamily\bfseries}             
{}                      
{0.5em}                 
{#1 #2{\hspace{0.2cm}--\hspace{-0.2cm}}  }   

\theoremstyle{mystyle}

\newtheorem{thm}{Theorem}
\newtheorem*{thm*}{Theorem}

\newtheorem{cor}[thm]{\hspace{-0.15cm}  {Corollary} }
\newtheorem{lem}[thm]{\hspace{-0.14cm}  {Lemma} }
\newtheorem{prop}[thm]{\hspace{-0.13cm} {Proposition}}
\newtheorem{defn}[thm]{ \hspace{-0.31cm} {Definition}}
\newtheorem*{rem*}{\hspace{-0.15cm} {Remark}}

\newtheoremstyle{mystyle2}
{3pt}               
{3pt}               
{\it }                      
{}                      
{\sffamily\bfseries}             
{}                      
{0.5em}                 
{\llap{#2 }#1{\hspace{0.2cm}--}}

\theoremstyle{mystyle2}

\newtheorem*{definition*}{Definition}
\newtheorem*{theorem*}{Theorem}
\newtheorem*{Remark*}{Remark}
\newtheorem*{lem*} {Lemma}
\newtheorem*{defn*} {Definition}
\newtheorem*{prop*} {Proposition}
\newtheorem*{cor*} {Corollary}

\newcommand{\ssk}{\smallskip}

\renewcommand{\epsilon}{\varepsilon}

\def\Labe{\mathfrak{e}}

\def\Labn{\mathfrak{n}}

\newcommand\bbN{\mathbb{N}}
\newcommand\bbR{\mathbb{R}}
\newcommand{\bbT}{\mathbb{T}}

\newcommand{\mcB}{\mathcal{B}} 
\newcommand{\mcC}{\mathcal{C}}

\newcommand{\mcH}{\mathcal{H}}
\newcommand{\mcI}{\mathcal{I}}

\newcommand{\mcM}{\mathcal{M}}

\newcommand\mcT{\mathcal T}

\begin{document}

\begin{center}
{\LARGE\sffamily{Parametrization of renormalized models for singular stochastic PDEs  \vspace{0.5cm}}}
\end{center}

\begin{center}
{\sf I. BAILLEUL}\footnote{I.B. acknowledges support from the CNRS/PIMS and UBC and the ANR-16-CE40-0020-01 grant.} \& {\sf Y. BRUNED}
\end{center}

\vspace{1cm}

\begin{center}
\begin{minipage}{0.8\textwidth}
\renewcommand\baselinestretch{0.7} \scriptsize \textbf{\textsf{\noindent Abstract.}} Let $\mathscr{T}$ be the regularity structure associated with a given system of singular stochastic PDEs. The paracontrolled representation of the $\sf \Pi$ map provides a linear parametrization of the nonlinear space of admissible models $\sf M=(g,\Pi)$ on $\mathscr{T}$, in terms of the family of para-remainders used in the representation. We give an explicit description of the action of the most general class of renormalization schemes on the parametrization space of the space of admissible models. The action is particularly simple for renormalization schemes associated with degree preserving preparation maps. The BHZ renormalization scheme has that property.
\end{minipage}
\end{center}

\vspace{0.6cm}


\section{Introduction}
\label{Introduction}

Regularity structures were introduced by M. Hairer \cite{Hai14} as a setting where to make sense and prove well-posedness of a large family of stochastic partial differential equations (PDEs) that come as scaling limit of microscopic discrete random dynamics where nonlinear and random effects balance each other. Each equation of this class is called a subcritical singular stochastic PDE. Prominent examples of subcritical singular stochastic PDE are given by the $1$-dimensional (KPZ) equation
$$
(\partial_t - \partial_x^2) u = \vert\partial_xu\vert^2 + \zeta,
$$
with a $(1+1)$-dimensional spacetime white noise $\zeta$, by the $2$ or $3$ dimensional parabolic Anderson model equation
$$
(\partial_t - \Delta) u = u\xi,
$$
with $\xi$ a space white noise, and by the $3$-dimensional scalar $(\Phi^4_3)$ equation from quantum field theory
$$
(\partial_t - \Delta)u = -u^3 + \zeta,
$$
with a $(1+3)$-dimensional spacetime white noise $\zeta$. Besides the fundamental works by M. Hairer and his co-authors \cite{Hai14, BHZ, ChandraHairer, BCCH18} a number of works provide partial views on different parts of the theory \cite{HairerBrazil, ChandraWeber, FrizHairer, HairerTakagi, CorwinShen}. We refer the reader to Bailleul \& Hoshino's {\it Tourist guide to regularity structures and singular stochastic PDEs} \cite{RSGuide} for a short self-contained reference on the algebraic and analytic sides of regularity structures theory and its applications to the study of singular stochastic PDEs.

One of the main features of the theory of regularity structures is the tight intertwining between its analytic, algebraic and probabilistic sides. To each subcritical singular stochastic PDE is attached an algebraic structure over which analytical objects are defined. Their realization as distributions on the state space of the equation requires the probabilistic construction of a finite number of distributions, a {\it model}. This construction involves an explicit construction with strong algebraic features, called {\it renormalization}. We tackle in the present work a question that is exactly at the intersection of the three sides of the subject: {\it Study the action on the parametrization space of the set of admissible models of the most general renormalization schemes}. To better grasp the stakes of that problem recall that the setting of regularity structures disentangles the task of solving an equation from the problem of making sense of a number of ill-defined quantities that are characteristic from the singular nature of the equation. The latter are encapsulated in the notion of model over a regularity structure. It provides a finite family of reference distributions/functions which are used to give local descriptions of possible solutions to a given singular stochastic PDE around each point in its state space. The construction of models associated with low regularity random noises is what renormalization is about. The systematic approach to the renormalization problem for singular stochastic partial differential equations (PDEs) was built gradually from Hairer's ad hoc construction in his groundbreaking work \cite{Hai14} to Bruned, Hairer and Zambotti's general setting for the BPHZ-type robust renormalization procedure \cite{BHZ} implemented by Chandra \& Hairer in \cite{ChandraHairer}. The dual action of this renormalization procedure on the equation was unveiled in Bruned, Chandra, Chevyrev and Hairer's work \cite{BCCH18}. The specific BHZ renormalization scheme was included in \cite{BrunedRecursive} by Bruned in a larger class of renormalization schemes, and the dual action of schemes of this class on the equation was investigated in Bailleul \& Bruned's work \cite{BailleulBruned} using algebraic insights from Bruned \& Manchon's work \cite{BrunedManchon}. 

The definition of a model $\sf (g,\Pi)$ over a given regularity structure $\mathscr{T}$ involves nonlinear operations that turn the metric space of models into a nonlinear space. Bailleul \& Hoshino were able in \cite{BH1, BH2} to provide a parametrization of the space of models over a given regularity structure by a linear space, a product of H\"older spaces. This parametrization involves the tools of paracontrolled calculus. Having such a parametrization is useful for understanding the structure of the space of models and \cite{BH1, BH2} contains a number of applications. The present work tackles the question of understanding the action of the most general renormalization scheme on the parametrization space of the models used for the study of systems of singular stochastic PDEs. The particular case of branched rough paths was investigated earlier by Tapia \& Zambotti in \cite{TZ18}  -- branched rough paths are a particular example of  models over a particular regularity structure, indexed by a time interval. Tapia \& Zambotti obtained a free transitive action of a product of H\"older spaces on the space of branched rough paths. The action of a general renormalization map on their parametrization space was investigated by Bruned in \cite{BrunedTZ}. However the particular case of branched rough paths only captures part of the structure of the general case.

\ssk

The regularity structures used for the study of singular stochastic PDEs have a particular structure described in depth in \cite{BHZ}. The models `adapted' to this structure are called {\it admissible}. We need a piece of notation to describe the parametrization of the set of admissible models over a given regularity structure $\mathscr{T}=\big((\mcT,\Delta),(\mcT^+,\Delta^+)\big)$. Given $\tau\in \mcT$ write
\begin{equation} \label{EqDelta}
\Delta\tau = \sum_{\sigma\leq\tau} \sigma\otimes(\tau/\sigma) \,\in \mcT\otimes \mcT^+.
\end{equation}
A choice of linear basis $\mcB$ of $\mcT$ fixes uniquely this decomposition by requiring that the elements $\sigma\in \mcT$ that appear in the sum belong to $\mcB$. This notation is only used in that sense in this work. In order to stick strictly to the statements proved in \cite{BH1} we formulate things in the case where the state space of the dynamics is the isotropic space $\bbR^d$ or its periodic version $\bbT^d$; this corresponds to elliptic equations. A similar result holds in the anisotropic setting used for the study of parabolic equations. The bilinear operator $\sf P$ below stands for a paraproduct operator; its definition or analytic properties are not needed in the present work, so we refer the reader to the first section of \cite{BH1} for more information. Let $\mathscr{T}$ be the BHZ regularity structure associated with a given (elliptic) singular stochastic PDE and let $\mcB$ be a basis of $\mcT$ -- details are given in Section \ref{SectionModelsPreparationMaps}. The following statement is a particular case of Theorem 2 in \cite{BH1}.

\medskip

\begin{thm} \label{ThmBH1}
Given any family of distributions $\big([\tau]\in\mcC^{\textsf{\emph{deg}}(\tau)}(\bbR^d)\big)_{\tau\in\mcB,  \textsf{\emph{deg}}(\tau)\leq 0}$, there exists a unique admissible model $\sf M = (g,\Pi)$ on $\mathscr{T}$ such that one has
$$
{\sf \Pi}\tau = \sum_{\sigma\leq\tau} {\sf P}_{{\sf g}(\tau/\sigma)}[\sigma],
$$
for all $\tau\in\mcB$ with $\textsf{\emph{deg}}(\tau) \leq 0$.
\end{thm}

\medskip

Note the specific form of the above representation of ${\sf \Pi}\tau$; a different paracontrolled representation of ${\sf \Pi}$ involving other functions than the ${\sf g}(\tau/\sigma)$ has for instance no a priori reason to give rise to a {\it parametrization} of the model. It is convenient to talk of a {\it bracket map} $[\,\cdot\,]$ associated with the model $\sf \Pi$. The precise statement of our main result involves notations that will be introduced below. We state it here in a qualitative form and refer the reader to Theorem \ref{ThmActionParametrizationBHZ} and Theorem \ref{ThmMain} for the full statements. The (degree preserving) preparation maps and their associated renormalization maps mentioned in Theorem \ref{ThmMain0} are defined in Definition \ref{DefnPreparationMap} and Definition \ref{DefnDegreePreserving} and equations \eqref{EqMWithI}, \eqref{EqMFromR} in Section \ref{SectionDegreePreserving}.

\medskip

\begin{thm} \label{ThmMain0}
Assume that an admissible model on $\mathscr{T}$ is given and parametrized by a family of distributions $[\tau]\in C^{\textsf{\emph{deg}}(\tau)}(\bbR^d)$, for $\tau\in\mcB$ with $\textsf{\emph{deg}}(\tau) \leq  0$. Let $R : \mcT\mapsto\mcT$ be a preparation map with associated renormalization map $M_{\!R} : \mcT\mapsto\mcT$ and renormalized model $\big({\sf g}^{\!R},{\sf \Pi}^{\!R}\big)$.   \vspace{0.15cm}

\begin{itemize}
   \item[\textbf{\textsf{(1)}}] If $R$ is degree preserving then the map ${\sf \Pi}^{\!R}$, hence the entire admissible model, is parametrized by the $[M_{\!R}\tau]$, for $\tau\in\mcB$ with $ \textsf{\emph{deg}}(\tau) \leq 0$.   \vspace{0.15cm}

   \item[\textbf{\textsf{(2)}}] In the general case of a non-degree preserving preparation map $R$ the bracket map $[\,\cdot\,]^{\!R}$ giving the parametrization of the map ${\sf \Pi}^{\!R}$ is given explicitly in terms of the bracket map $[\,\cdot\,]$.
   \end{itemize}   
\end{thm}

\medskip

Item \textbf{\textsf{(1)}} means that the renormalized model $({\sf g}^R,{\sf \Pi}^R)$ is characterized by the fact that one has
$$
{\sf \Pi}^{\!R}\tau = \sum_{\sigma\leq\tau} {\sf P}_{{\sf g}^{\!R}(\tau/\sigma)}[M_R\sigma]
$$
for all $\tau\in\mcB$ with $\textsf{\emph{deg}}(\tau) \leq 0$. The description of the renormalized model is not as nice in the general setting of item \textbf{\textsf{(2)}}. We emphasize here that the class of degree preserving preparation maps is much larger than the class of BHZ renormalization maps. We start Section \ref{SectionModelsPreparationMaps} by giving back the main features of the BHZ regularity structure associated with a given system of singular stochastic PDEs. The renormalization schemes that we consider in this work are built from maps called preparation maps. A number of useful results about these maps and their associated renormalization maps are given in Section \ref{SectionModelsPreparationMaps}. The proof of Theorem \ref{ThmMain0} in the particular case of degree preserving preparation maps is the object of Section \ref{SectionDegreePreserving}; the general case is treated in Section \ref{SectionGeneral}.

\medskip

\section{Basics on regularity structures and  preparation maps}
\label{SectionModelsPreparationMaps}

We first recall in Section \ref{SubsectionBHZ} the setting of BHZ regularity structures that we use in the present work. Preparation maps and their elementary properties are described in Section \ref{SubsectionPreparationMaps}.

\medskip

\subsection{BHZ regularity structures associated with singular stochastic PDEs}
\label{SubsectionBHZ}

Consider for simplicity the case of a single equation whose mild formulation writes
\begin{equation} \label{mild}
u = K \ast \Big(F(u,\nabla u, \ldots) \xi + G(u,\nabla u, \ldots)\Big),
\end{equation}
with $\ast$ standing for space or spacetime convolution and $\xi$ denoting a non-constant `noise'. (What follows works verbatim when working with several equations and more noises.) The main idea of regularity structures is to iterate the mild formulation locally when the noises $\xi$ is replaced by a regularized version. Nonlinearities $F$ or $G$ are Taylor-expanded around arbitrary base points $x \in \mathbb{R}^d$ and one can construct recentered iterated integrals around each point. One obtains a local description of the solution of \eqref{mild} of the form
\begin{equation} \label{B-series}
u(y) = \sum_{\tau \in T} c_{\tau}(x) \left( {\sf \Pi}_x \tau \right)(y) + R(x,y)
\end{equation}
where the collection $T$ consists of combinatorial objects called decorated trees that we will be described below. The coefficients $c_{\tau}(x)$ are Taylor-type coefficients, the $ {\sf \Pi}_x \tau $ are iterated integrals recentered around the point $ x $ and $ R(x,y) $ is a Taylor-type remainder. 

\ssk

We now introduce decorated trees and their symbolic notations. Pick two symbols $\mcI$ and $\Xi$ and let $ \mathcal{D} := \lbrace \mcI,\Xi \rbrace \times \mathbb{N}^d$ define the set of edge decorations. Decorated trees over $ \mathcal{D} $ are triples of the form  $\tau_{\Labe}^{\Labn} =  (\tau,\Labn,\Labe) $ where $\tau$ is a non-planar rooted tree with node set $N_\tau$ and edge set $E_\tau$. The maps $\Labn : N_\tau \rightarrow \mathbb{N}^d$ and $\Labe : E_\tau \rightarrow \mathcal{D}$ are node, respectively edge, decorations. The set of decorated trees is denoted by $T_0$ and we write $\mathcal{T}_0$ for the linear span of $T_0$. The index $0$ refers to the fact that $T_0$ and $ \mathcal{T}_0$ are background objects from which the spaces $\mcT$ and $\mcT^+$, part of the regularity structure associated with \eqref{mild}, will be defined. The tree product $ \cdot $ on $ T_0$ is defined by 
\begin{equation}  \label{treeproduct}
 	(\tau,\Labn,\Labe) \cdot (\tau',\Labn',\Labe') 
 	= (\tau \cdot \tau',\Labn + \Labn', \Labe + \Labe')\;, 
\end{equation} 
where $\tau \cdot \tau'$ is the rooted tree obtained by identifying the roots of $ \tau$ and $\tau'$. The sums $ \Labn + \Labn'$ and  $\Labe + \Labe'$ mean that decorations are added at the root and extended to the disjoint union by setting them to vanish on the other tree. Each edge and vertex of both trees keeps its decoration, except the roots which merge into a new root decorated by the sum of the previous two decorations. 

\ssk

{\color{gray} $\bullet$} We will use mainly in this work a symbolic notation for these decorated trees. Denote by $\{e_1, \ldots, e_d\}$ the canonical basis of $ \mathbb{N}^d$.

\begin{enumerate}
   \item[--] An edge decorated by  $ (\mcI,a) \in \mathcal{D} $  is denoted by $ \mathcal{I}_{a} $. The symbol $  \mathcal{I}_{a} $ is also viewed as  the operation that grafts a tree onto a new root via a new edge with edge decoration $ a $. The new root at hand remains decorated with 0.  The operator $ \mathcal{I}_a $ encodes a space-time convolution with the kernel $\partial^a K $.
   
   \item[--]  An edge decorated by $ (\Xi,0) \in \mathcal{D} $ is denoted by $ \Xi $. We only consider decorated trees having such edge as a terminal edge meaning that one of its extremity is a leaf.  We also disregard trees having some decorations of the form $(\Xi,a)$ with $a \neq 0$.

   \item[--] A factor $ X^k $ encodes a single node  $ \bullet^{k} $ decorated by $ k \in \mathbb{N}^d$. We write $ X_i $, $ i \in \lbrace 1,\ldots,d\rbrace $, to denote $ X^{e_i} $. The element $ X^0 $ is identified with the empty tree $ \textsf{\textbf{1}} $.
 \end{enumerate}
Any decorated tree $ \tau $ admits the following decomposition 
 \[
	\tau =  X^{k} \Xi^m \prod_{i=1}^{n} \mathcal{I}_{a_i}(\tau_i) ,
 \]
 where the $\tau_i $ are decorated trees and the product $ \prod_{i=1}^n $ is the tree product.   The factor $ X^{k} $ expresses the fact that the root of $\tau$ is decorated by $k$.  Using symbolic notation one can reformulate the tree product \eqref{treeproduct} as
 \begin{equation*}
	\left( X^{k} \Xi^m \prod_{i} \mathcal{I}_{a_i}(\tau_i)  \right) 
	\left( X^{ k'} \Xi^{m'} \prod_{j} \mathcal{I}_{b_j}(\sigma_j) \right) 
	= X^{k + k'} \Xi^{m+m'}\prod_{i} \mathcal{I}_{a_i}( \tau_i) \prod_{j}  \mathcal{I}_{b_j}(\sigma_j).
\end{equation*}
We note here for later use that a forest is a collection of trees equipped with the forest product given by the disjoint union. Recall that characters are linear maps that are multiplicative.

\ssk

{\color{gray} $\bullet$} We now associate to decorated trees numbers that depend on their decorations. These numbers form a degree map denoted by $\textsf{deg} : T_0 \rightarrow \mathbb{R}$ defined inductively by the following relations
\begin{equation} \label{deg_def}
\textsf{deg}(\textbf{\textsf{1}}) = 0, \quad \textsf{deg}(\Xi) = \alpha, \quad \textsf{deg}(\sigma) = \textsf{deg}(\tau) + \beta - |a|, \quad 
 \textsf{deg}(\sigma\tau) = \textsf{deg}(\sigma) + \textsf{deg}(\tau),
\end{equation} 
where $ \alpha < 0 $ is the (space or spacetime) regularity of the noise $ \Xi $ in a suitable Hölder space, $ \beta $ corresponds to the gain of regularity in Schauder estimate for the convolution operator with the kernel $K$. We also denote by $ |\cdot|_{\Xi} : T_0 \rightarrow \mathbb{N} $ the map that counts the number of noises in any given decorated tree. Using the degree map $\textsf{deg}$ we associate to any  $E \subset T_0$ the set
\begin{equation} \label{Tplus}
	 E^+ 
 	:=   \Big\lbrace X^{k} \prod_{i=1}^{n} \mathcal{I}^{+}_{a_i}(\tau_i)\,;\,
	\textsf{deg}(\mathcal{I}_{a_i}(\tau_i)) > 0, \tau_i \in E, \, k \in \mathbb{N}^d \Big\rbrace . 
\end{equation}
This definition means that all the branches outgoing from the root must be of positive degree. We use a different symbol $ \mathcal{I}_a^+ $ to stress that $ E^+ $ is not a subset of $ E $ as there is no constraint on the edges connected to the root for elements of $ E^+ $. We define $\mathcal{T}_0^+$ as the linear span of $T_0^+$ and equip $ \mathcal{T}_0^+ $ with a Hopf algebra structure. Its product is the tree product and  its coproduct $ \Delta^{\!+} $  is given by 
\begin{equation} \label{EqDelta+}
\Delta^{\!+}(\mcI^{+}_a\tau) := \sum_{\ell \in \mathbb{N}^d} \bigg(\mcI^+_{a+\ell}\otimes\frac{(-X)^{\ell}}{\ell !}\bigg)\Delta \tau + \textsf{\textbf{1}}\otimes  \mcI^{+}_a\tau,
\end{equation}
and its antipode map $ S^{+} $ is given inductively by the relation
\begin{equation} \label{EqS+}
S^+ (\mcI^{+}_a\tau) = - \sum_{\ell \in \mathbb{N}^d}  \mathcal{M}^+ \left( \mcI^{+}_{a+\ell} \otimes  \frac{X^{\ell}}{\ell!}S^{+}\right) \Delta \tau.
\end{equation}
Denote by $\mcM^+ : \mcT_0^+\otimes \mcT_0^+\rightarrow \mcT_0^+$ the multiplication operator in $\mcT_0^+$ and  by $ \textsf{\textbf{1}}^{\star} : \mathcal{T}_0 \rightarrow \mathbb{R} $ the counit linear map equal to one on $ \textsf{\textbf{1}} $ and zero otherwise. The main identities that we will use in the sequel are the co-associativity and the characterization of the antipode given below by
\begin{equation} \label{coa_deltap}
\left( \Delta^{\!+} \otimes \textrm{Id} \right) \Delta^{\!+} = \left( \textrm{Id} \otimes \Delta^{\!+}   \right) \Delta^{\!+}
\end{equation}
and
\begin{equation} \label{char_antipode}
\mathcal{M}^+ \left( S^{+} \otimes \textrm{Id} \right) \Delta^{\!+} = 
\mathcal{M}^+ \left(  \textrm{Id} \otimes S^{+}  \right) \Delta^{\!+} =  \textsf{\textbf{1}}^{\star} \textsf{\textbf{1}}
\end{equation}
The space $ \mathcal{T}_0$ is equipped with a co-action $ \Delta $ defined by 
\begin{equation} \label{coaction} \begin{split} 
\Delta(\bullet) &:= \bullet\otimes \textsf{\textbf{1}}, \quad \textrm{ for }\bullet\in\big\{\textsf{\textbf{1}}, X_i, \Xi\big\}, \quad \Delta X_i = X_i \otimes \textsf{\textbf{1}} + \textsf{\textbf{1}} \otimes X_i,   \\
\Delta(\mcI_a\tau) &:= (\mcI_a\otimes\textrm{Id})\Delta \tau + \sum_{\vert\ell +m\vert < \textsf{deg}(\mcI_a\tau)} \frac{X^\ell}{\ell!}\otimes \frac{X^m}{m!}\mcI^+_{a+\ell+m}(\tau).
\end{split}\end{equation}
This definition turns the pair $(\mcT_0,\Delta)$ into a right comodule over $\mathcal{T}_0^+$ and one has the following compatibility identity
\begin{equation} \label{coaction_identity}
\left( \Delta \otimes \textrm{Id} \right) \Delta = \left( \textrm{Id} \otimes \Delta^{\!+}   \right) \Delta.
\end{equation}
(These expressions used in Hairer's original work \cite{Hai14} are different from the expressions used by Bruned, Hairer and Zambotti in \cite{BHZ}. One moves from \cite{Hai14} to \cite{BHZ} by performing a change of basis in $\mcT_0^+$ and taking $\sum_{\ell \in \mathbb{N}^{d}}  \frac{(-X)^{\ell}}{\ell !} \mathcal{I}^+_{a + \ell}(\tau)$ in the role of $\mcI^+_a(\tau)$. The induction rule giving the action of $\Delta$ on the abstract integration operator is more useful here in the form of relation \eqref{EqDefnDeltaMCirc} than in the form given in \cite{BHZ}, identity (3.6) in \cite{RSGuide}.)

\ssk

{\color{gray} $\bullet$} The BHZ regularity structure associated to equation \eqref{mild} is a subset $ T \subset T_0$ defined from the product appearing from the right hand side of the equation. (The acronym `BHZ' is chosen after the names of the three authors of \cite{BHZ}.) For example, due to the affine structure of the noise, the term $ \Xi^2 $ does not appear in any meaningful formal expansion of a potential solution to equation \eqref{mild}; this puts a constraint on the way the decorated trees for this equation are constructed. Such constraints are formalized through the notion of {\it normal complete rule}. We  refer the reader to Section 5 of \cite{BHZ} where those rules have been detailed. The main property of $\mathcal{T}_0$ is that the co-module and Hopf algebraic structures satisfied by $\mcT_0$ and $ \mcT_0^+$ are also satisfied by $\mcT$ (linear span of $ T $) and $\mcT^+$ (linear span of $ T^+ $) with the {\it same maps} $\Delta$ and $\Delta^{\!+}$. Suitable assumptions on the products of \eqref{mild}, called {\it local sub-criticality}, guarantee that $ \mathcal{T} $ admits a direct sum decomposition $\mcT=\bigoplus_{\beta\in A}\mcT_\beta$ involving finite dimensional vector spaces $\mcT_\beta$ generated by decorated trees of degree $\beta$. We write $\mathscr{T}$ for the pair $\big((\mcT,\Delta),(\mcT^+,\Delta^+)\big)$ that defines the (BHZ) regularity structure associated with equation \eqref{mild}. It is clear from this description that $\mcT$ and $\mcT^+$ come equipped with canonical bases.

\ssk

{\color{gray} $\bullet$} We recall that the notion of admissibility of a model $\sf (g,\Pi)$ over the BHZ regularity structure for equation \eqref{mild} is relative to the operator $K$ and that admissible models satisfy 
$$
{\sf \Pi}(\mcI_a\tau) = (D^aK)*({\sf \Pi}\tau)
$$
and 
\begin{equation} \label{EqRelationAdmissibleModel}
{\sf g}_x^{-1}\big(\mathcal{I}^+_a\tau\big) = - \big(D^{a} K * {\sf \Pi}_x\tau\big)(x),
\end{equation}
where 
\begin{equation} \label{def_model}
 {\sf \Pi}_x  := \left( {\sf \Pi} \otimes {\sf g}^{-1}_x \right) \Delta,
\end{equation}
for all $x$ in the state space. 

\ssk

{\color{gray} $\bullet$} The algebraic structure given by a regularity structure encodes the mechanics of local expansions and re-expansions for the functions involved in the analysis of a given equation. The renormalization procedure involved in the definition of a model associated with a low regularity noise is encoded in another algebraic structure. Given a subset $E$ of $ T_0$, we denote by $E^-$ the forest formed of elements in $E$ having negative degree. One uses a coproduct $ \Delta^{\!-} $ and a co-action $\delta : \mcT\rightarrow\mcT^-\otimes\mcT$ to construct renormalization maps parametrized by the group of characters $ \mathcal{G}^- $ of $ \mathcal{T}^- $. The group structure on the set of characters is derived from the fact that $\mathcal{T}^-$ has a Hopf algebra structure when equipped with $\Delta^{\!}$ and an appropriate antipode map $S^-$. The map $\delta$ turns the space $\mcT$ into a left comodule over $\mathcal{T}^-$. Here is the key identity in this business
\begin{equation*} 
\left( \textrm{Id} \otimes \delta \right) \delta =  \left( \Delta^{\!-}\otimes \textrm{Id}  \right) \delta,
\end{equation*}
It gives the formula
\begin{equation*}
\ell_1 \star \ell_2 := \left( \ell_1 \otimes \ell_2\right) \Delta^{\!-},\qquad \ell^{-1} = \ell(S^- \cdot),
\end{equation*}
for the convolution product $ \star $ of two characters of $\mcT^-$ and their inverse, and provides an action of such characters on $\mcT$
\begin{equation} \label{EqMEll}
M_{\ell} := \left( \ell \otimes  \textrm{Id} \right) \delta, \quad M_{\ell} \circ M_{g} = M_{\ell \star g}.
\end{equation}
Adding extra decorations one can obtain a simple action of the maps $M_{\ell}$ on admissible models $\sf (g,\Pi)$  setting
\begin{equation} \label{simple_action}
{\sf \Pi}_x^{M_{\ell}} := {\sf \Pi}_x M_{\ell}.
\end{equation}
(The action of $M_\ell$ on the $\sf g$-part of the model will be described later.) Such a definition is possible due to the co-interaction between $ \Delta  $ and $ \delta $ described by
\begin{equation*} \label{cointeraction}
\mathcal{M}^{(13)(2)(4)} \left( \delta \otimes \delta \right) \Delta = \left( \textrm{Id} \otimes \Delta \right)  \delta
\end{equation*}
where $ \delta $ is aslo defined as a map from $ \mathcal{T}^+ $ into $ \mathcal{T}^{-} \otimes \mathcal{T}^+ $ and $\mathcal{M}^{(13)(2)(4)}$ is given for $ \tau_1, \tau_3 \in \mathcal{T}^{-}  $, $ \tau_2 \in \mathcal{T} $ and $ \tau_4 \in \mathcal{T}^+ $ by:
\begin{equation*}
\mathcal{M}^{(13)(2)(4)} \left(  \tau_1 \otimes \tau_2 \otimes \tau_3 \otimes \tau_4 \right)
= \tau_1 \tau_3 \otimes \tau_2 \otimes \tau_4
\end{equation*}
where the product between $ \tau_1 $ and $ \tau_3 $ is the forest product. We will continue in the sequel with a more general formalism for the renormalization procedure introduced by Bruned in \cite{BrunedRecursive}.

\medskip

We will denote below by $\mcB$ the canonical basis of $\mcT$; this fixes in particular the shorthand notation in formula \eqref{EqDelta} describing $\Delta$.

\medskip

\subsection{Preparation maps}
\label{SubsectionPreparationMaps}

We take from \cite{BrunedRecursive} the following definition.

\medskip

\begin{defn} \label{DefnPreparationMap}
A \textbf{\textsf{preparation map}} is a map 
$$
R : \mcT\rightarrow \mcT
$$ 
that fixes polynomials and such that 
\begin{itemize}
   \item[--] for each $ \tau \in T $ there exist finitely many $\tau_i \in T$ and constants $\lambda_i$ such that
\begin{equation} \label{EqAnalytical}
R \tau = \tau + \sum_i \lambda_i \tau_i, \quad\textrm{with}\quad \textsf{\emph{deg}}(\tau_i) \geq \textsf{\emph{deg}}(\tau) \quad\textrm{and}\quad |\tau_i|_{\Xi} < |\tau|_{\Xi},
\end{equation} 

   \item[--] one has 
 \begin{equation} \label{EqCommutationRDelta}
 ( R \otimes \textrm{\emph{Id}}) \Delta = \Delta R.
 \end{equation}
\end{itemize}
\end{defn}

\medskip

Identity \eqref{EqAnalytical} gives an upper triangular structure to preparation maps that is useful for inductive proofs. Identity \eqref{EqCommutationRDelta} encodes a commutation property between the recentering operator encoded in the map $\Delta$ and the `renormalization' operator encoded in the map $R$.

\medskip

\noindent \textbf{\textsf{Example --}} {\sl The archetype of a preparation map is defined from a map $\delta_r$, with the index `$r$' for `root', defined similarly as the splitting map $\delta$, but extracting from any $\tau\in T$ only one diverging subtree of $\tau$ with the same root as $\tau$ at a time, and summing over all possible such subtrees -- see Definition 4.2 in \cite{BrunedRecursive}. Given a character $\ell$ of the algebra $\mcT^-$ the map
\begin{equation} \label{EqREll}
R_\ell := (\ell\otimes\textrm{\emph{Id}})\delta_r
\end{equation}
is a preparation map.   }   \hfill $\RHD$ 

\medskip

We will work exclusively with preparation maps $R : \mcT\rightarrow \mcT$ such that 
$$
R\,\mcI_a=\mcI_a,
$$ 
for all $a$. Let $M_{\!R}^{\!\times} : \mcT\rightarrow \mcT$ and $M_{\!R} : \mcT\rightarrow \mcT$ be the maps uniquely defined from $R$ by requiring that $M_{\!R}^{\!\times}$ is {\it multiplicative} and satisfies
\begin{equation} \label{EqMWithI}
M_{\!R}^{\!\times}(\mcI_a\tau) = \mcI_a\big(M_{\!R}^{\!\times} (R\tau)\big)
\end{equation}
and 
\begin{equation} \label{EqMFromR}
M_{\!R} := M^{\!\times}_{\!R} R.
\end{equation}
The map $M_{\!R}$ is the {\it renormalization map associated with the preparation map $R$}. While this map is {\it not} multiplicative it follows from \eqref{EqMFromR} that $M_{\!R}$ commutes with all the integration operators $\mcI_a$. Note that the map $M_{R_\ell}$ associated with \eqref{EqREll} is of the type introduced in \cite{BHZ}. In the setting of \cite{BHZ} the structure of the renormalization schemes on $\mcT$ and their induced actions on $\mcT^+$ are encoded in the splitting maps $\delta : \mcT\rightarrow\mcT^-\otimes\mcT$ and characters of the algebra $\mcT^-$. Here the algebraic structure associated with the renormalization map $M_{\!R}$ is entirely encoded in the latter. The following result is used to describe the renormalized model; it was first proved in Proposition 8.36 in \cite{Hai14}. We give an elementary proof to be self-contained.

\medskip

\begin{lem} \label{LemInvertibility}
The map
$$
\big(\textrm{\emph{Id}}\otimes\mcM^+\big)\,(\Delta\otimes\textrm{\emph{Id}}) : \mcT\otimes \mcT^+\rightarrow \mcT\otimes \mcT^+
$$ 
is invertible.
\end{lem}

\medskip

\begin{Dem}
Writing 
$$
\Delta\sigma = \sum_{\sigma_1\leq\sigma} \sigma_1\otimes \sigma/\sigma_1,
$$
one has for $\sigma\in \mcT$ and $\tau\in \mcT^+$
$$
\big(\textrm{Id}\otimes\mcM^+\big)\,(\Delta\otimes\textrm{Id}) (\sigma \otimes\tau) = \sum_{\sigma_1\leq\sigma} \sigma_1\otimes \big(\sigma/\sigma_1\tau\big)
$$
and the only element in the previous sum whose $\mcT^+$-component has maximum degree is $\sigma\otimes\tau$. This shows the injectivity of the map $\big(\textrm{Id} \otimes\mcM^+\big)\,(\Delta\otimes \textrm{Id})$. It surjectivity comes from the fact that 
$$
\big(\textrm{Id}\otimes\mcM^+\big)\,(\Delta\otimes\textrm{Id})(\sigma\otimes\tau) =: \sigma\otimes\tau + N(\sigma\otimes\tau),
$$ 
for a nilpotent map $N$, so a Neumann series gives the inverse of $(\textrm{Id}\otimes\mcM^+)\,(\Delta\otimes\textrm{Id})$. A different representation 
\begin{equation} \label{EqInverse}
\big(\textrm{Id}+N\big)^{-1} = (\textrm{Id}\otimes\mcM^+)\big(\textrm{Id}\otimes S^+\otimes\textrm{Id}\big)\big(\Delta\otimes\textrm{Id}\big)
\end{equation}
was proved by Bruned in Lemma 3.20 of \cite{BrunedRecursive}. Formula \eqref{EqInverse} plays a role in the proof of Lemma \ref{LemPreserving} in Section \ref{Section::3}.
\end{Dem}

\medskip

It follows from Lemma \ref{LemInvertibility} that one defines inductively two maps
$$
\delta_{\!R} : \mcT\rightarrow \mcT\otimes \mcT^+,  \quad  M_{\!R}^+ : \mcT^+\rightarrow \mcT^+,
$$
setting
\begin{equation} \label{EqConditionDeltaM}
(\textrm{Id}\otimes\mcM^+)(\Delta\otimes\textrm{Id})\delta_{\!R} := (M_{\!R}\otimes M_{\!R}^{\!+})\Delta,
\end{equation}
with $M_{\!R}^{\!+} : \mcT^+\rightarrow \mcT^+$, the {\it multiplicative} map fixing the monomials and such that one has
$$
M_{\!R}^{\!+}\Big(\mcI^+_a(\tau)\Big) = \mcM^+\big(\mcI^+_a\otimes \textrm{Id}\big)\delta_{\!R}\tau,
$$
for all $\tau\in \mcT$. Lemma \ref{LemFactorizationDeltaR} below gives a useful representation of the map $\delta_R$ that leads to a direct proof of Proposition \ref{LemmaTriangular}. The proof of Lemma \ref{LemFactorizationDeltaR} is the very place where we take advantage of the fact that we work with renormalization maps built from a preparation map, as opposed to working with a general renormalization map as those of Section 8.3 of Hairer' seminal work \cite{Hai14}. Define a {\it multiplicative} map 
$$
\delta_{\!R}^{\!\times} : \mcT\rightarrow \mcT\otimes\mcT^+
$$ 
setting
\begin{equation} \label{EqDefnDeltaMCirc} \begin{split} 
\delta_{\!R}^{\!\times}(\bullet) &:= \bullet\otimes \textsf{\textbf{1}}, \quad \textrm{ for }\bullet\in\big\{\textsf{\textbf{1}}, X_i, \Xi\big\},   \\
\delta_{\!R}^{\!\times}(\mcI_a\tau) &:= (\mcI_a\otimes\textrm{Id})\delta_{\!R}^{\!\times}(R\tau) - \sum_{\vert\ell\vert\geq \textsf{deg}(\mcI_a\tau)} \frac{X^\ell}{\ell!}\otimes\mcM^+\big(\mcI^+_{a+\ell}\otimes\textrm{Id}\big)\delta_{\!R}^{\!\times}(R\tau).
\end{split}\end{equation}
 
\medskip

\begin{lem} \label{LemFactorizationDeltaR}
One has $\delta_{\!R} = \delta_{\!R}^{\!\times} \, R$.
\end{lem}

\medskip

\begin{Dem} We proceed by induction. Using identity \eqref{EqMFromR} to write 
\begin{equation*}
(M_{\!R}\otimes M_{\!R}^+)\Delta = (M_{\!R}^{\!\times}\otimes M_{\!R}^+)\Delta R
\end{equation*}
and the fact that $ R $ is invertible we are down to checking that one has 
\begin{equation*}
(\textrm{Id}\otimes\mcM^+)(\Delta\otimes\textrm{Id})\delta_{\!R}^{\!\times} = (M_{\!R}^{\!\times}\otimes M_{\!R}^+)\Delta.
\end{equation*}
It suffices by multiplicativity to consider a tree of the form $\mathcal{I}_{a}(\tau)$, for which one has on the one hand
\begin{equation*}
(M_{\!R}^{\!\times}\otimes M_{\!R}^{\!+})\Delta \mathcal{I}_a(\tau) \overset{\eqref{coaction}}{=} (\mcI_a\otimes\textrm{Id})\big(M_{\!R}^{\!\times}\otimes M_{\!R}^{\!+}\big)\Delta \tau + \sum_{\vert\ell +m\vert < \textsf{deg}(\mcI_a\tau)} \frac{X^\ell}{\ell!}\otimes \frac{X^m}{m!} \, M_{\!R}^{\!+}\big(\mcI^+_{a+\ell+m}(\tau)\big).
\end{equation*}
On the other hand we have
\begin{equation*}
\begin{aligned}
(\textrm{Id}\otimes\mcM^+)(\Delta\otimes\textrm{Id})\delta_{\!R}^{\!\times}\big(\mathcal{I}_a(\tau)\big) &\overset{\eqref{EqDefnDeltaMCirc}}{=}
(\mcI_a\otimes\textrm{Id})(\textrm{Id}\otimes\mcM^+)(\Delta\otimes\textrm{Id})\delta_{\!R} \tau   \\ 
&\quad + \sum_{\ell, m \in \mathbb{N}^d} \frac{X^\ell}{\ell!}\otimes \frac{X^m}{m!} \mcM^+\big(\mcI^+_{a+\ell+m}\otimes\textrm{Id}\big)\delta_{\!R}\tau   \\
&\quad- \sum_{\vert\ell + m\vert\geq \textsf{deg}(\mcI_a\tau)} \frac{X^\ell}{\ell!}\otimes\frac{X^m}{m!} \mcM^+\big(\mcI^+_{a+\ell+m}\otimes\textrm{Id}\big)\delta_{\!R} \tau
\end{aligned}
\end{equation*}
We conclude by applying the induction hypothesis on $ \tau $. A similar proof was performed in Proposition 3.19 of \cite{BrunedRecursive} using the explicit formula \eqref{EqInverse} for $\big(\textrm{Id}+N\big)^{-1}$.

\end{Dem}

\medskip

\begin{defn*}
A map $A : \mcT\rightarrow \mcT\times \mcT^+$, with $A\tau = \sum\tau_1\otimes\tau_2$, is said to be \textbf{\textsf{upper triangular}} if $\textsf{\emph{deg}}(\tau_1)\geq \textsf{\emph{deg}}(\tau)$, for all $\tau_1$ in the preceding decomposition of $A\tau$, $ \tau \in T $.  
\end{defn*}

\medskip

\begin{prop} \label{LemmaTriangular}
The map $\delta_{\!R}$ is upper triangular.
\end{prop}

\medskip

\begin{Dem}
It suffices from the property \eqref{EqAnalytical} of preparation maps to see that $\delta_{\!R}^{\!\times}$ is upper triangular. This point is obtained from \eqref{EqDefnDeltaMCirc} and the multiplicativity of $\delta_{\!R}^{\!\times}$ by an elementary induction on $\textsf{deg}(\tau)+\vert\tau\vert_{\Xi}$.
\end{Dem}

\medskip

It follows from Proposition \ref{LemmaTriangular} and the definition of $M_{\!R}^{\!+}$ that $\textsf{deg}\big(M_{\!R}^{\!+}\sigma\big)\geq \textsf{deg}(\sigma)$, for all $\sigma\in T^+$. The last ineguality means that for $ M_{\!R}^{\!+}\sigma = \sum_i \lambda_i \sigma_i $, with $\textsf{deg}\big(\sigma_i \big)\geq \textsf{deg}(\sigma)$ for every $i$.

\bigskip

\section{The case of degree preserving preparation maps}
\label{SectionDegreePreserving}
\label{Section::3}

We prove the first part of Theorem \ref{ThmMain0} in this section. Degree preserving preparation maps are defined below in Definition \ref{DefnDegreePreserving}. The algebraic properties  enjoyed by the renormalization maps associated with the class of degree preserving preparation maps allow a direct construction of renormalized admissible models close to what is done for the BHZ models from Bruned-Hairer-Zambotti's work \cite{BHZ}. The construction involved in the general case is not as simple; it will be detailed in Section \ref{SectionGeneral}.

\medskip

Recall from Theorem \ref{ThmBH1} that the formula 
\begin{equation} \label{EqPCRepresentation}
{\sf \Pi}\tau = \sum_{\sigma\leq \tau} {\sf P}_{{\sf g}(\tau/\sigma)}[\sigma],
\end{equation}
for $\tau\in T$ with $\textsf{deg}(\tau)\leq 0$, provides a parametrization of the set of admissible models over a large class of regularity structures containing those used for the study of singular stochastic PDEs -- the BHZ regularity structures from Section \ref{SubsectionBHZ}. Let us stress that if we are given an admissible model $\sf (g,\Pi)$, formula \eqref{EqPCRepresentation} defines uniquely the bracket map $[\cdot]$. Indeed identity
\begin{equation}
{\sf \Pi}\tau = \sum_{\sigma < \tau} {\sf P}_{{\sf g}(\tau/\sigma)}[\sigma] + [\tau],
\end{equation} 
shows that $[\tau]$ depends on $ {\sf \Pi} \tau $ and the $[\sigma]$ and $ {\sf g} $ applied to elements which are strictly smaller. The renormalization maps used in \cite{BHZ} are built from specific features of BHZ regularity structures  and from a character on the Hopf algebra $(\mcT^-,\Delta^{\!-})$ that is in co-interaction with $(\mcT,\Delta)$, such as encoded in \eqref{EqMEll} and \eqref{simple_action}. A single feature of the fine structures involved in the definition of the preparation map associated with the BHZ renormalization map is of importance here. It singles out a large class of preparation maps for which the action of their associated renormalization maps on the parametrization space takes the simple form given in Theorem \ref{ThmActionParametrizationBHZ} below. The BHZ renormalization maps form one family of this class.

\medskip

\begin{defn} \label{DefnDegreePreserving}
A preparation map is said to be \textbf{\textsf{degree preserving}} if for each $\tau\in T$ there exists finitely many $\tau_i\in T$ and constants $\lambda_i$ such that
\begin{equation} \label{EqAnalytical_2}
R \tau = \tau + \sum_i \lambda_i \tau_i, \quad\textrm{with}\quad \textsf{\emph{deg}}(\tau_i) = \textsf{\emph{deg}}(\tau) \quad\textrm{and}\quad |\tau_i|_{\Xi} < |\tau|_{\Xi}.
\end{equation}
\end{defn}

\medskip

Compare condition \eqref{EqAnalytical_2} with condition \eqref{EqAnalytical} involved in the definition of an arbitrary preparation map. The introduction in \cite{BHZ} of decorated trees with extended decorations allows precisely to design a setting where the splitting map associated with the renormalization procedure enjoys a similar property. (One works in this setting with two degree maps $\textsf{deg}$ and $\textsf{deg}_-$, with $\textsf{deg}_-$ not taking into account the extended decorations and involved in the definition of the Hopf algebra $(\mcT^-,\Delta^-)$ that is part of the renormalization structure on $\mathscr{T}$.) Although elementary it is of fundamental importance that the maps $M_{\!R}^{\!\times}$ and $M_{\!R}$ associated to a degree preserving preparation map $R$ are also degree preserving. This is what allows to prove the next statement by induction.

\medskip

\begin{lem} \label{LemPreserving}
For any degree preserving preparation map $R$ one has
\begin{equation} \label{EqSimpleDeltaR}
\delta_{\!R}\tau=(M_{\!R}\tau)\otimes \textsf{\textbf{1}}
\end{equation}
and the co-interaction identity 
\begin{equation} \label{cointeraction}
\Delta M_{\!R} = \big(M_{\!R}\otimes M_{\!R}^{\!+}\big)\Delta.
\end{equation}
One further has that $M_{\!R}^{\!+}$ commutes with the antipode $S^{\!+}$.
\end{lem}

\medskip

\begin{Dem}
Note that the co-interaction identity \eqref{cointeraction} is equivalent from \eqref{EqCommutationRDelta} to the identity
\begin{equation} \label{EqCointeraction}
\Delta M^{\!\times}_{\!R}  = \big(M_{\!R}^{\!\times}\otimes M_{\!R}^{\!+}\big)\Delta.
\end{equation}
This identity involves only multiplicative maps on $\mcT$, so it suffices to prove it for elements of $\mcT$ of the form $\Xi, X^k$ or $\mcI_a(\tau)$. It is elementary to check it for $\Xi$ and $X^k$. We prove identities \eqref{EqSimpleDeltaR} and \eqref{EqCointeraction} for elements of $\mcT$ of the form $\mcI_a(\tau)$ by induction on $\textsf{deg}(\tau)+\vert\tau\vert_\Xi$, for a generic $\tau\in T$. We use the symbol $(\star)$ above an $=$ sign to emphasize the use of the induction assumption in a sequence of equalities. Write 
$$
M_{\!R}\tau=\tau+\sum_i c_i \sigma_i,
$$ 
for constants $c_i$, with $\vert\sigma_i\vert_\Xi<\vert\tau\vert_\Xi$. As
$$
\textsf{deg}(\mcI_a\tau) > \textsf{deg}(\tau),
$$ 
for all $\mcI_a\tau\in T$, one has
$$
\textsf{deg}(\mcI_a\tau)+\vert\tau\vert_\Xi > \textsf{deg}(\sigma_i)+\vert \sigma_i\vert_\Xi, \qquad\forall\,i
$$
from the fact that $M_{\!R}$ is degree preserving. This justifies the use of the induction hypothesis in the $(\star)$ equality below.
\begin{equation*} \begin{aligned}
\Delta M^{\!\times}_{\!R}(\mathcal{I}_a\tau)  & = \Delta  \mathcal{I}_a(M_{\!R} \tau) =   (\mcI_a\otimes\textrm{Id})\Delta (M_{\!R} \tau) + \sum_{\vert\ell +m\vert < \textsf{deg}(\mcI_a\tau)} \frac{X^\ell}{\ell!}\otimes \frac{X^m}{m!}\mcI^+_{a+\ell+m}(M_{\!R} \tau)
\\ &\hspace{-0.05cm}\overset{(\star)}{=}  (\mcI_a M_{\!R} \otimes M_{\!R}^+)\Delta  \tau+ \sum_{\vert\ell +m\vert < \textsf{deg}(\mcI_a\tau)} \frac{X^\ell}{\ell!}\otimes \frac{X^m}{m!}M^+_{\!R}\mcI^+_{a+\ell+m}( \tau)
\\ & =  \big(M_{\!R}^{\!\times}\otimes M_{\!R}^{\!+}\big)\Delta \mathcal{I}_a(\tau)
\end{aligned} \end{equation*}
The bound on $ |\ell + m| $ in the first line comes from the degree preserving property of $ M_R $: one has  $  \textsf{deg}(\mcI_a\tau) =  \textsf{deg}(\mcI_a\sigma_i)$  for all $ i $.
We have used the induction assumption about \eqref{EqCointeraction} for the first term in the right hand side of the third equality and the induction assumption about \eqref{EqSimpleDeltaR} for the second term in the right hand side of that equality coupled with the fact that
$$
\mcI^+_b( M_{\!R} \tau) = \mcM^+\big(\mcI^+_b\otimes \textrm{Id}\big)\delta_{\!R}\tau = M_{\!R}^{\!+}\big(\mcI^+_b\tau\big), \qquad \forall\,\tau\in T.
$$
Identity \eqref{EqConditionDeltaM} defining $\delta_{\!R}$ then reads 
$$
(\textrm{Id}\otimes\mcM^+)(\Delta\otimes\textrm{Id})\delta_{\!R} \sigma = (M_{\!R}\otimes M_{\!R}^{\!+})\Delta \sigma = \Delta (M_{\!R} \sigma),
$$
and it follows from the explicit formula \eqref{EqInverse} that 

\begin{equation*}
\begin{aligned}
\delta_{\!R} \sigma & = (\textrm{Id}\otimes\mcM^+)\big(\textrm{Id}\otimes S^+\otimes\textrm{Id}\big)\big(\Delta\otimes\textrm{Id}\big) \Delta(M_{\!R} \sigma)
\\ & = (\textrm{Id}\otimes\mcM^+)\big(\textrm{Id}\otimes S^+\otimes\textrm{Id}\big)\big( \textrm{Id} \otimes\Delta^{\!+}\big) \Delta (M_{\!R} \sigma)
\\ & = (  \textrm{Id} \otimes \textsf{\textbf{1}}^{*} \textsf{\textbf{1}})  \Delta (M_{\!R} \sigma) =  M_{\!R} \sigma \otimes \textsf{\textbf{1}}.
\end{aligned}
\end{equation*}
where we have used the following property of the antipode $ S^+ $ 
\begin{equation*}
\mathcal{M}^+\left( S^{+} \otimes \textrm{Id} \right) \Delta^{\!+}= \textsf{\textbf{1}}^{*} \textsf{\textbf{1}}
\end{equation*}
One sees that $ M^{\!+}_{\!R}$ and $S^+$ commute using the inductive relation 
$$
S^+\mcI_a(\tau) = - \sum_{\ell\in\bbN^d} \mcM^+\big(\mcI^+_{a+\ell}\otimes \frac{X^\ell}{\ell!}\,S^+\big)\Delta \tau
$$
for the antipode $S^+$ and writing

\begin{equation*} \begin{aligned}
 M^{\!+}_{\!R}  S^+ (\mcI^{+}_a\tau) & = - \sum_{\ell \in \mathbb{N}^d}  \mathcal{M}^+ \left(  M^{\!+}_{\!R}\mcI^{+}_{a+\ell} \otimes  \frac{X^{\ell}}{\ell!} M^{\!+}_{\!R} S^{+}\right) \Delta \tau = - \sum_{\ell \in \mathbb{N}^d} \mathcal{M}^+ \left(  \mcI^{+}_{a+\ell} M_R \otimes  \frac{X^{\ell}}{\ell!} S^{+} M^{\!+}_{\!R}\right) \Delta \tau
 \\ & = - \sum_{\ell \in \mathbb{N}^d}   \mathcal{M}^+ \left(\mcI^{+}_{a+\ell} \otimes  \frac{X^{\ell}}{\ell!} S^{+} \right) \Delta M_{\!R} \tau = S^{+}  \big(\mcI^{+}_a M_{\!R} \tau\big) = S^{+} M^{\!+}_{\!R}  (\mcI^{+}_a \tau).
 \end{aligned} \end{equation*}
\end{Dem}

\medskip

Similar computations are involved in Remark 4.2.6 and Proposition 4.2.8 of \cite{bruned:tel-01306427}. Note that it follows from \eqref{EqSimpleDeltaR} that the multiplicative map $M_{\!R}^{\!+}$ satisfies in that case the relation
$$
M_{\!R}^{\!+}\big(\mathcal{I}^+_a(\tau)\big) =   \mathcal{I}^+_a( M_{\!R} \tau).
$$
Since $M_{\!R}$ is degree preserving we read on the previous identity that $M_{\!R}^{\!+}$ is also degree preserving. One then proves similarly as in the proof of Lemma \ref{LemPreserving} that $M^{\!+}_{\!R}$ satisfies the co-interaction identity
\begin{equation} \label{EqCoInteractionPlus}
\big(M_{\!R}^{\!+}\otimes M_{\!R}^{\!+}\big)\Delta^{\!+} = \Delta^{\!+} M^{\!+}_{\!R}.
\end{equation}
Given an admissible model $\sf (g,\Pi)$ on $\mathscr{T}$ set
$$
{\sf g}^{\!R} := {\sf g}\circ M_{\!R}^{\!+}, \quad {\sf \Pi}^{\!R} := {\sf \Pi}\circ M_{\!R}.
$$
It follows from the fact that $M_{\!R}^{\!+}$ commutes with the antipode $S^+$ that 
$$
({\sf g}^{\!R})^{-1} = {\sf g}^{-1}\circ M_{\!R}^{\!+}.
$$

\medskip

\begin{cor}
The pair $\big({\sf g}^{\!R}, {\sf \Pi}^{\!R}\big)$ defines an admissible model on $\mathscr{T}$.
\end{cor}

\medskip

\begin{Dem}
On the one hand, identities \eqref{EqConditionDeltaM} and \eqref{EqSimpleDeltaR} ensure that
\begin{equation} \label{EqPiMx} \begin{split}
{\sf \Pi}^{\!R}_x\tau &= \big({\sf \Pi}^{\!R}\otimes ({\sf g}^{\!R}_x)^{-1}\big)\Delta\tau = \Big\{{\sf \Pi}^{\!R}\otimes\big({\sf g}_x^{-1}\circ M_{\!R}^{\!+}\big)\Big\}\Delta\tau = \big({\sf \Pi}\otimes {\sf g}_x^{-1}\big)(M_{\!R}\otimes M_{\!R}^{\!+})\Delta\tau   \\
			 &\hspace{-0.16cm}\overset{\eqref{EqConditionDeltaM}}{=} \big({\sf \Pi}_x\otimes {\sf g}_x^{-1}\big)\,\delta_{\!R}\tau \overset{\eqref{EqSimpleDeltaR}}{=} {\sf \Pi}_x(M_{\!R}\tau).
\end{split} \end{equation}
It follows from this identity and the fact that $M_{\!R}$ is degree preserving that ${\sf \Pi}^R_x$ satisfies the analytic estimates required from a model on $\mathscr{T}$. On the other hand, the co-interaction identity \eqref{EqCoInteractionPlus} gives
\begin{equation*} \begin{aligned}
{\sf g}^{\!R}_{yx} &= \Big({\sf g}^{\!R}_y\otimes({\sf g}^{\!R}_x)^{-1}\Big)\Delta^{\!+} = \big({\sf g}_y\otimes{\sf g}_x^{-1}\big)\Big(M_{\!R}^{\!+}\otimes M_{\!R}^{\!+}\Big)\Delta^{\!+}  
		\\ &\hspace{-0.14cm}\overset{\eqref{EqCoInteractionPlus}}{=} \big({\sf g}_y\otimes{\sf g}_x^{-1}\big)\Delta^{\!+} M^{\!+}_{\!R} =  {\sf g}_{yx} \circ M_{\!R}^{\!+}.
\end{aligned} \end{equation*}
It follows from this identity and the fact that $M_{\!R}^{\!+}$ is degree preserving that ${\sf g}^{\!R}_{yx}$ satisfies the analytic estimates required from a model on $\mathscr{T}$.
\end{Dem}

\medskip

\begin{thm} \label{ThmActionParametrizationBHZ}
Assume $\mathscr{T}=\big((\mcT,\Delta),(\mcT^+,\Delta^{\!+}),(\mcT^-,\Delta^{\!-})\big)$ is the BHZ regularity structure associated with a system of singular stochastic PDEs. Let $\mcB$ stands for the canonical linear basis of $\mcT$ and let $\sf (g,\Pi)$ be an admissible model on $\mathscr{T}$, with associated bracket map $[\,\cdot\,]$ in its paracontrolled representation \eqref{EqPCRepresentation}. For any degree preserving preparation map $R$ the admissible model $\big({\sf g}^{\!R}, {\sf \Pi}^{\!R}\big)$ on $\mathscr{T}$ is parametrized by the family $\big( [M_{\!R}\tau]\in C^{\textsf{\emph{deg}}(\tau)}(\bbR^d)\big)_{\tau\in \mcB,\, \textsf{\emph{deg}}(\tau) \leq 0}$. 
\end{thm}

\medskip

\begin{Dem}
The action of $M_{\!R}$ on the parametrization set of the space of admissible models is given by
$$
{\sf \Pi}^{\!R}\tau = {\sf \Pi}(M_{\!R}\tau) \overset{\eqref{cointeraction}}{=} \sum_{\textsf{\textbf{1}} <\sigma\leq\tau} {\sf P}_{{\sf g}(M_{\!R}^+(\tau/\sigma))}[M_{\!R}\sigma] = \sum_{\textsf{\textbf{1}} <\sigma\leq\tau} {\sf P}_{{\sf g}^{\!R}(\tau/\sigma)}[M_{\!R}\sigma].
$$
The second equality follows from the fact that $M_{\!R}^{\!+}$ commutes with the antipode $S^{\!+}$ and from the formula \eqref{EqgR} for ${\sf g}^{\!R}$. The term $ {\sf P}_{f} 1 $ is equal to zero for any $ f \in \mathcal{S}'(\mathbb{R}^d) $. Therefore one can remove the term $ \sigma = \textsf{\textbf{1}}$ in the sum giving ${\sf \Pi}^{\!R}\tau$. The fact that $M_{\!R}\sigma$ is a sum of terms of the same degree as the degree of $\sigma$ shows that the preceding identity gives a parametrization of the model associated with ${\sf \Pi}^{\!R}$ by the $[M_{\!R}\sigma]$, for all $\sigma$ with negative degree.
\end{Dem}

\medskip

Theorem~\ref{ThmActionParametrizationBHZ} provides a nice action of $ M_R $ on the parametrisation space. Indeed one gets the explicit expression
\begin{equation} \label{explicit_brackets}
[\cdot]^{M_R}  = [M_R \cdot]
\end{equation}
for the renormalised brackets $ [\cdot]^{M_R}$. This is an important result because it shows that the recursive definition of the $ [\cdot] $ is decouple from the renormalization and it is connected with \eqref{simple_action}. In both cases, the co-interaction \eqref{cointeraction} between the  recentering given by $  \Delta $ and the renormalization given by $ \delta $  plays a major role in the proof. In the context of degree preserving renormalization maps, \eqref{cointeraction} is replaced by \eqref{EqCointeraction}. In the next section, we will consider a more general set up and one gets a weaker result saying that $ [\cdot]^{M_R} $ depends recursively on the brackets $ [\cdot] $.

\medskip

\noindent \textbf{\textsf{The example of branched rough paths --}} {\sl An action of a renormalization group was observed previously in Bruned's work \cite{BrunedTZ} on the renormalization of branched rough paths. This kind of rough paths was introduced by Gubinelli in \cite{GubinelliBranched}. Hairer \& Kelly showed in \cite{HairerKelly} that they can be seen as weak geometric rough paths over a larger space. See e.g. Cass \& Weidner's work \cite{CassWeidner} or Bailleul's work \cite{BailleulCassWeidner} for a quick grasp on branched rough paths.

The regularity of a branched rough path is quantified by an exponent $\gamma\in(0,1)$, and a $\gamma$-branched rough path is indexed by decorated trees $\tau$; denote by $|\tau|$ the number of nodes in $\tau$. Fix $\gamma\in(0,1)$, and for a continuous function $h$ on $[0,1]$ write $h_t$ for its value at time $t$. Tapia and Zambotti exhibited in \cite{TZ18} a free transitive action of the product space of H\"older spaces
$$
\mcH^\gamma := \bigg\{g=\big(g(\tau)\big)_{\tau\in\mcB}\in\prod_{\tau\in\mcB, |\tau|\leq 1/\gamma}C^{\gamma\vert\tau\vert}([0,1])\,;\, g_0(\tau) = 0 \bigg\},
$$
where $\mcB$ is a certain collection of $\gamma$-branched rough paths, on the space of all $ \gamma $-branched rough paths. One of the main results of Bruned's work \cite{BrunedTZ} provides an explicit formula for the map $g^M\in\mcH^\gamma$ sending any $\gamma$-branched rough path $X$ to $X \circ M$, for a renormalization map $M$ associated in this particular setting to a preparation map of BHZ type, hence degree preserving map. The map $g^M$ is given in Theorem 4.4 of \cite{BrunedTZ} and takes the form
\begin{equation} \label{formula11}
g^M_t(\tau) - g^M_s(\tau) = \big\langle \overline{ X_{ts} \circ M},\tau \big\rangle - \langle \overline{X}_{ts},\tau\rangle.
\end{equation}
where $ \overline{X}$ is the Lyons-Victoir extension of $X$. The latter is not so explicit, so Theorem \ref{ThmActionParametrizationBHZ} above gives a better description of the action of a renormalization map even in that setting. The paracontrolled parametrization bypasses in particular the problem emphasized in Remark 4.6 of \cite{BrunedTZ} related to the nonlinear character of the Lyons-Victoir extension map. Bellingeri, Friz, Paycha \& Preiss' recent work \cite{BFPP} contains material related to the question of renormalization of smooth rough paths.   }

\bigskip

\section{The general case}
\label{SectionGeneral}

We prove the second part of Theorem \ref{ThmMain0} in this section. In the particular case of degree preserving preparation maps Lemma \ref{LemPreserving} gives a simple form for $\delta_{\!R}$, one has the co-interaction identity \eqref{cointeraction} and the commutation of $M_{\!R}^{\!+}$ with the antipode $S^+$. These properties do not hold in the case of a general preparation map $R$ so one cannot use the mechanics of the proof of Theorem \ref{ThmActionParametrizationBHZ} in the general case. One can however give an explicit description of the admissible model ${\sf M}^{\!R} = \big({\sf g}^{\!R}, {\sf \Pi}^{\!R}\big)$ associated with $R$ and infer from it an inductive description of the bracket map $[\,\cdot\,]^{\!R}$ associated with ${\sf \Pi}^{\!R}$. We describe the admissible model associated with a preparation map in Section \ref{SubsectionConstructionModel} before describing the bracket map $[\,\cdot\,]^{\!R}$ in Section \ref{SubsectionGeneraResult}.

\medskip

\subsection{Renormalised model associated with a preparation map}
\label{SubsectionConstructionModel}

We will use in the next statement a density argument in the space of models that requires the introduction of a regularity structure $\mathscr{T}^{(\epsilon)}$, indexed by a positive regularity exponent $\epsilon$. The only difference between $\mathscr{T}^{(\epsilon)}$ and $\mathscr{T}$ is the notion of degree $\textsf{deg}^{(\epsilon)}$ on $\mathscr{T}^{(\epsilon)}$, defined as $\textsf{deg}^{(\epsilon)}(\tau) = \textsf{deg}(\tau)-\epsilon$, for all $\tau\in\mcT$ or $\tau\in\mcT^+\backslash\{\textbf{\textsf{1}} \}$. The exponent $\epsilon$ is chosen small enough $\textsf{deg}^{(\epsilon)}(\tau)$ to be positive for all $\mcT^+\backslash\{\textbf{\textsf{1}}\}$. Given now an admissible model $\big({\sf g}, {\sf \Pi}\big)$ on $\mathscr{T}$, set for all $\tau\in \mcT$ and $\sigma\in \mcT^+$
\begin{equation} \label{EqDefnModelM}
{\sf \Pi}^{\!R}\tau := {\sf \Pi}\big(M_{\!R}\tau\big), \quad ({\sf g}^{\!R})^{-1}(\sigma) := {\sf g}^{-1}\big(M_{\!R}^{\!+}\sigma\big).
\end{equation}
The map ${\sf \Pi}^{\!R}$ satisfies the admissibility condition from the fact that $M_{\!R}$ commutes with the operators $\mcI_a$ and from the admissibility of the map ${\sf \Pi}$.

\medskip

\begin{prop} \label{PropRenormalizedModel}
The pair $\big({\sf g}^{\!R},{\sf \Pi}^{\!R}\big)$ defines an admissible model on $\mathscr{T}^{(\epsilon)}$.
\end{prop}

\medskip

We have in particular 
\begin{equation} \label{EqgR}
{\sf g}^{\!R}(\sigma) = {\sf g}\big(S^+M_{\!R}^+S^+\sigma\big), \quad\forall\,\sigma\in\mcT^+.
\end{equation}
We will see as a corollary of Theorem \ref{ThmMain} that $\big({\sf g}^{\!R},{\sf \Pi}^{\!R}\big)$ is actually a model on $\mathscr{T}$. Bruned has proved in Section 3 of \cite{BrunedRecursive} a version of Proposition \ref{PropRenormalizedModel} for {\it continuous admissible models}. The use below of a density argument allows to extend the result to {\it all admissible models}.

\medskip

\begin{Dem}
Smooth models are models for which all the $\sf \Pi\tau$ and ${\sf g}(\sigma)$ are smooth functions. We know from Theorem 2 in \cite{BH1} or Theorem 5 in \cite{BH2}, giving paracontrolled parametrization of the space of admissible models, that the set of smooth admissible models on $\mathscr{T}$ is dense in the topology associated with the canonical injection of the set of models on $\mathscr{T}$ in the set of models on $\mathscr{T}^{(\epsilon)}$. See also Theorem 2.14 in Singh and Teichmann's work \cite{SinghTeichmann} for a similar statement. It suffices then to prove that for any smooth model $\sf (g,\Pi)$ the pair $\big({\sf g}^{\!R},{\sf \Pi}^{\!R}\big)$ defines an admissible model on $\mathscr{T}$ -- this is what we prove in the following.

\ssk

Identity \eqref{EqConditionDeltaM} ensures that
\begin{equation} \label{EqPiMx} \begin{split}
{\sf \Pi}^{\!R}_x\tau &= \big({\sf \Pi}^{\!R}\otimes ({\sf g}^{\!R}_x)^{-1}\big)\Delta\tau = \Big\{{\sf \Pi}^{\!R}\otimes\big({\sf g}_x^{-1}\circ M_{\!R}^{\!+}\big)\Big\}\Delta\tau   \\
		   	 &= \big({\sf \Pi}\otimes {\sf g}_x^{-1}\big)(M_{\!R}\otimes M_{\!R}^{\!+})\Delta\tau   \\
 		   	 &\hspace{-0.16cm}\overset{\eqref{EqConditionDeltaM}}{=}\big({\sf \Pi}_x\otimes {\sf g}_x^{-1}\big)\,\delta_{\!R}\tau.
\end{split} \end{equation}
It follows from this identity and the fact that $\delta_{\!R}$ is upper triangular, Lemma \ref{LemmaTriangular}, that ${\sf \Pi}^R_x$ satisfies the analytic estimates required from a model on $\mathscr{T}$. Note that this holds for {\it all} admissible models $\sf (g,\Pi)$, smooth or not. (This point will be used in the proof of Corollary \ref{CorModelT}.)

\ssk

Note that it follows from \eqref{EqPiMx} and the admissibility of the model $\sf (g,\Pi)$ that one has for all $\tau\in\mcT$ and all $x$
\begin{equation} \label{EqFormulagxRMoins1} \begin{aligned}
\big({\sf g}^{\!R}_x\big)^{-1}\big(\mcI^+_a\tau\big) & = {\sf g}_x^{-1}\Big(M_{\!R}^{\!+}\big(\mcI^+_a\tau\big)\Big)   \\ 
&= {\sf g}_x^{-1}\Big(\mcM^+\big(\mcI^+_a\otimes \textrm{Id}\big)\delta_{\!R}\tau\Big)   \\ 
&\hspace{-0.16cm}\overset{\eqref{EqRelationAdmissibleModel}}{=} \Big((- D^aK * {\sf \Pi}_x)(x)\otimes {\sf g}_x^{-1}\Big)\delta_{\!R}\tau   \\ 
&\hspace{-0.22cm}\overset{\eqref{EqPiMx}}{=} - \Big(D^aK * {\sf \Pi}^{R}_x \tau\Big)(x).
 \end{aligned} \end{equation}
Write for all $x,y$
\begin{equation*}
{\sf g}^{\!R}_{yx} := \Big({\sf g}^{\!R}_x\otimes \big({\sf g}_y^{\!R}\big)^{-1}\Big) \Delta^{\!+},
\end{equation*}
and define now a multiplicative map from $\mcT$ into itself setting for all $\tau\in T$
$$
\widehat{{\sf g}^{\!R}_{yx}}(\tau) := \big(\textrm{Id}\otimes {\sf g}^{\!R}_{yx}\big)\Delta\tau.
$$
Denoting by $\mu_\beta$ the component of any $\mu\in \mcT$ in $\mcT_\beta$ in the grading $\bigoplus_{\beta\in A}\mcT_\beta$ of $\mcT$, the analytic estimates required from ${\sf g}^{\!R}_{yx}$ for $\big({\sf g}^{\!R},{\sf\Pi}^{\!R}\big)$ to be a model on $\mathscr{T}$ are equivalent to having
\begin{equation} \label{EqEstimategyx}
\Big\vert\big(\widehat{{\sf g}^{\!R}_{yx}}(\sigma)\big)_\beta\Big\vert \lesssim \vert y-x\vert^{\textsf{deg}(\sigma)-\beta}
\end{equation}
for all $\sigma\in T$ and all $\beta$ with $\beta<\textsf{deg}(\tau)$, for all $x,y$. We have
$$
\widehat{{\sf g}^{\!R}_{yx}}(X_i) = X_i + (x_i-y_i){\bf 1},\quad \widehat{{\sf g}^{\!R}_{yx}}(\Xi) = \Xi.
$$
The following identity is where working with smooth models helps -- continuous models would make the job as well.

\ssk

\begin{lem}
One has the identity
\begin{equation} \label{recursivePi}
\widehat{{\sf g}^{\!R}_{yx}}\big(\mcI_a\tau\big) = \mcI_a\Big(\widehat{{\sf g}^{\!R}_{yx}}\tau\Big) - \sum_{\vert\ell\vert< \textsf{\emph{deg}}(\mcI_a\tau)} \frac{(X+x-y)^\ell}{\ell !}\,{\sf \Pi}_x^{\!R}\Big(\mcI_{a+\ell}\big(\widehat{{\sf g}^{\!R}_{yx}}\tau\big)\Big)(y).
\end{equation}
\end{lem}

\ssk

\noindent \textbf{\textsf{Proof  --}} Note the pointwise evaluation of ${\sf \Pi}_x^{\!R}$ at a given point $y$; we work with smooth models to make sense of it -- having a continuous model would be sufficient. We briefly recall how one can obtain \eqref{recursivePi} as the settings in \cite{BH1} and \cite{BrunedRecursive} are not striclty speaking the same. The inductive relation \eqref{coaction} on $\Delta$ gives
\begin{equation} \label{EqFirstFormula} \begin{aligned}
\widehat{{\sf g}^{\!R}_{yx}}\big(\mcI_a\tau\big) &=  \mcI_a\Big(\widehat{{\sf g}^{\!R}_{yx}}\tau\Big)+ \sum_{\vert\ell + m\vert< \textsf{deg}( \mcI^+_a\tau)} \frac{X^\ell}{\ell !}\, {\sf g}^{\!R}_{yx}\bigg(\frac{X^m}{m!}\,\mcI^+_{a+ \ell + m}(\tau)\bigg)   \\ 
&= \mcI_a\Big(\widehat{{\sf g}^{\!R}_{yx}}\tau\Big) +\sum_{\vert k\vert< \textsf{deg}( \mcI^+_a\tau )} \frac{(X+x-y)^k}{k!}\, {\sf g}^{\!R}_{yx}\big(\mcI^+_{a+ k}(\tau)\big)
\end{aligned} \end{equation}
Rewriting ${\sf g}^{\!R}_{yx}$ under the form
\begin{equation*}
{\sf g}^{\!R}_{yx} =  \Big(\big(({\sf g}^{\!R}_x)^{-1} \circ S^+\big) \otimes ({\sf g}_y^{\!R})^{-1} \Big) \Delta^{\!+},
\end{equation*}
in order to use relation \eqref{EqFormulagxRMoins1} giving $({\sf g}^{\!R}_x)^{-1}$ and $({\sf g}^{\!R}_y)^{-1}$, one has

\begin{equation*} \begin{aligned}
& {\sf g}^{\!R}_{yx}\big(\mathcal{I}_b(\tau)\big) \\ &\hspace{-0.16cm}\overset{\eqref{EqS+}}{=} ({\sf g}^{\!R}_{y})^{-1}\big(\mathcal{I}^+_b(\tau)\big) - \sum_{m,n \in \mathbb{N}^d} \bigg( \Big\{({\sf g}^{\!R}_x)^{-1}\mcI^+_{b+m+n}\otimes \frac{(-x)^m}{m!} ({\sf g}^{\!R}_x)^{-1}\Big\} S^+ \Delta\otimes\frac{y^{n}}{n !} ({\sf g}_y^{\!R})^{-1} \bigg)\Delta \tau   \\
&= ({\sf g}^{\!R}_{y})^{-1}\big(\mathcal{I}^+_b(\tau)\big) - \sum_{k \in \mathbb{N}^d} \frac{(y-x)^{k}}{k!}\Big(({\sf g}^{\!R}_x)^{-1}\mcI^+_{b+k}\otimes \Big(({\sf g}^{\!R}_x)^{-1} S^+\otimes ({\sf g}_y^{\!R})^{-1} \Big)\Delta^{\!+} \Big)\Delta \tau   \\ 
&= ({\sf g}^{\!R}_{y})^{-1}\big(\mathcal{I}^+_b(\tau)\big)  - \sum_{k \in \mathbb{N}^d} \frac{(y-x)^{k}}{k!} \Big(({\sf g}^{\!R}_x)^{-1}\mcI^+_{b+k}\otimes {\sf g}^{\!R}_{yx}  \Delta^{\!+} \Big)\Delta \tau
 \\ &  = ({\sf g}^{\!R}_{y})^{-1}\big(\mathcal{I}^+_b(\tau)\big) - \sum_{k \in \mathbb{N}^d}  \frac{(y-x)^{k}}{k !} \,({\sf g}^{\!R}_{x})^{-1}\Big(\mathcal{I}^+_{b + k}\big(\widehat{{\sf g}^{\!R}_{yx}}(\tau)\big)\Big),
\end{aligned} \end{equation*}
for all $b$. One gets identity \eqref{recursivePi} as follows from the preceding equality using the explicit expression for $({\sf g}_x^{\!R})^{-1}$ given in \eqref{EqFormulagxRMoins1} and the relation ${\sf \Pi}^{\!R}_x\circ \widehat{{\sf g}^{\!R}_{yx}} = {\sf \Pi}^{\!R}_y$. Write $ \widehat{{\sf g}^{\!R}_{yx}}( \tau) = \sum_i \lambda_i \tau_i $, with $ \textsf{deg}(\tau_i)  \leq \textsf{deg}(\tau) $ -- note that ${\sf \Pi}^R_{y} \tau = {\sf \Pi}^R_{x} \widehat{g^{\!R}_{yx}} \tau = \sum_i \lambda_i {\sf \Pi}^R_{x} \tau_i $, then
\begin{equation*} \begin{split}
{\sf g}^{\!R}_{yx}\big(\mcI^+_{a+ \ell}(\tau)\big) &= - \Big((D^{a+\ell} K)*({\sf \Pi}^{\!R}_y \tau)\Big)(y) - \sum_i \lambda_i \sum_{m \in \mathbb{N}^{d}} \frac{(y-x)^m}{m!} \big({\sf g}^{\!R}_x\big)^{-1}\big(\mathcal{I}^{+}_{a+\ell+m}(  \tau_i)\big)   \\
&= \sum_i \lambda_i  \left\{ - \Big((D^{a+\ell} K)*({\sf \Pi}^{\!R}_x \tau_i)\Big)(y) - \sum_{m \in \mathbb{N}^{d}} \frac{(y-x)^m}{m!} \big({\sf g}^{\!R}_{x}\big)^{-1}\big(\mathcal{I}^{+}_{a+\ell+m}(  \tau_i)\big)\right\} 
\end{split} \end{equation*}

\[
\begin{aligned}
&= \sum_i \lambda_i  {\sf \Pi}_x^{\!R}\big(\mathcal{I}^{+}_{a+\ell}(\tau_i)\big)(y) = {\sf \Pi}_x^{\!R}\Big(\mathcal{I}^{+}_{a+\ell}\big(\widehat{{\sf g}^{\!R}_{yx}}( \tau)\big)\Big)(y).
\end{aligned}
\] 
Together with identity \eqref{EqFirstFormula}, this gives the recursive formula \eqref{recursivePi}.   \hfill $\RHD$

\ssk

An elementary induction as in the proof of Proposition 3.16 in \cite{BrunedRecursive} shows then that the size bound \eqref{EqEstimategyx} on $\widehat{{\sf g}^{\!R}_{yx}}(\sigma)$ holds for all $\sigma=\mcI_a(\tau)$, with $\tau\in T$. This works as follows. Let $\alpha< \textsf{deg}( \sigma)$. If $ \alpha \in \mathbb{R} \setminus \mathbb{N} $, let us write $\widehat{{\sf g}^{\!R}_{yx}}\big( \tau \big)=\tau+\sum_i \lambda^i_{yx} \tau_i$ with 
 \[
\textsf{deg}(\tau_i) = \alpha_i<  \textsf{deg}(\tau), \qquad \big|(\tau_i)_{\alpha_i} \big| \lesssim | x-y|^{\textsf{deg}(\tau)-\alpha_i} ;
 \]
 then if 
\[
\big| \big(\widehat{{\sf g}^{\!R}_{yx}}(\sigma) \big)_{\alpha}\big| = \big| \big( \mathcal{I}_a\big(\widehat{{\sf g}^{\!R}_{yx}}(\tau)\big) \big)_{\alpha}\big| \lesssim\sum_i \mathds{1}_{\lbrace\alpha_i + \beta -  |a|= \alpha\rbrace} | x -y |^{\textsf{deg}(\tau) + \beta - |a|-\alpha_i} \lesssim | x-y|^{\textsf{deg}(\sigma)- \alpha}
\]
Now, if $ \alpha \in \mathbb{N} $ and $ \alpha < \textsf{deg}(\sigma)$ then 
\begin{equation*} \begin{split}
& \big|\big(\widehat{{\sf g}^{\!R}_{yx}}(\sigma)\big)_{\alpha}\big|   = \left\vert \left( 
 \sum_{ \alpha \leq \vert\ell\vert< \textsf{deg}(\sigma)} \frac{(X+x-y)^{\ell}}{\ell !} {\sf \Pi}^{\!R}_{x}\Big(\mathcal{I}_{a+\ell}\big( \widehat{{\sf g}^{\!R}_{yx}}(\tau)\Big)(y) \right)_{\alpha} \right\vert   \\
 & \lesssim \sum_{ \alpha \leq | \ell | <  \textsf{deg}(\sigma)}  |x-y|^{ |\ell|-\alpha} \sum_{\gamma \leq  \textsf{deg}(\tau )} |x-y|^{\gamma -|\ell| +\beta - |a| } \, \big| \big(\widehat{{\sf g}^{\!R}_{yx}}(\tau)\big)_{\gamma } \big|   \\
\end{split} \end{equation*} 
\begin{equation*} \begin{split} 
& \lesssim \sum_{ \alpha \leq | \ell | < \textsf{deg}(\sigma)}  |x-y|^{ |\ell|-\alpha} \sum_{\gamma \leq  \textsf{deg}( \tau )} |x-y|^{\gamma -|\ell| + \beta -| a | } |x-y |^{\textsf{deg}(\tau) - \gamma}  \lesssim | x - y |^{\textsf{deg}(\sigma) - \alpha}.
\end{split} \end{equation*} 
The multiplicativity of $\widehat{{\sf g}^{\!R}_{yx}}$ on $\mcT$ ensures that the bound \eqref{EqEstimategyx} holds for all $\sigma\in \mcT$.
\end{Dem}

\medskip

\begin{Remark*}
Proposition \ref{PropRenormalizedModel} can be proved in a different way, defining first a map $\delta_{\!R}^{\!+} : \mcT^+\rightarrow\mcT^+$ via the identity 
$$
(\textrm{\emph{Id}}\otimes\mcM^+)(\Delta^{\!+}\otimes\textrm{\emph{Id}})\delta_{\!R}^{\!+} = \big(S^+M_{\!R}^{\!+}S^+\otimes M_{\!R}^{\!+}\big)\Delta^+.
$$
One can prove as in Lemma \ref{LemInvertibility} that the map $(\textrm{\emph{Id}}\otimes\mcM^+)(\Delta^{\!+}\otimes\textrm{\emph{Id}}) : \mcT^+\otimes\mcT^+\rightarrow \mcT^+\otimes\mcT^+$ is invertible. The defining relation for $\delta_{\!R}^{\!+}$ ensures that 

\begin{equation*} \begin{split}
{\sf g}^{\!R}_{yx} &= \big({\sf g}^{\!R}_y\otimes({\sf g}^{\!R}_x)^{-1}\big)\Delta^{\!+}  = \big({\sf g}_y\otimes{\sf g}_x^{-1}\big)\Big(S^+M_{\!R}^{\!+}S^+\otimes M_{\!R}^{\!+}\Big)\Delta^{\!+}   \\
			 &= \big({\sf g}_{yx}\otimes ({\sf g}_x)^{-1}\big)\delta_{\!R}^{\!+}.
\end{split} \end{equation*}
A deep and fairly non-trivial result of Hairer \& Quastel ensures that the map $\delta_{\!R}^{\!+}$ is upper triangular if the map $\delta_{\!R}$ is upper triangular -- see Lemma B.1 in \cite{HQ}. The size estimates on ${\sf g}^{\!R}_{yx}(\sigma)$, for any $\sigma\in\mcT^+$, follows then from the preceding formula for ${\sf g}^{\!R}_{yx}$ and Hairer \& Quastel's result. It shows directly that $\big(g^{\!R},{\sf \Pi}^{\!R}\big)$ is an admissible model on $\mathscr{T}$ at the price of using Lemma B.1 of \cite{HQ} as a blackbox. Our proof is elementary and does not use Hairer \& Quastel's result; it follows the proof of Theorem 3.19 in \cite{BrunedRecursive}. We recover in Corollary \ref{CorModelT} below the fact that $\big(g^{\!R},{\sf \Pi}^{\!R}\big)$ is an admissible model on $\mathscr{T}$ rather than just an admissible model on $\mathscr{T}^{(\epsilon)}$.
\end{Remark*}

\bigskip

\subsection{Parametrization of renormalized models}
\label{SubsectionGeneraResult}

Assume now that we work with any preparation map $R$ on $\mcT$. The co-interaction identity \eqref{cointeraction} and the fact that $M_{\!R}^{\!+}$ commutes with $S^+$ are not guaranteed to hold anymore so the proof of Theorem \ref{ThmActionParametrizationBHZ} breaks down. We use instead the map $\delta_{\!R}$ which provides a connection between the renormalized model and the original model. We use the shorthand notation 
\begin{equation} \label{EqDescriptionDeltaM}
\delta_{\!R}\tau =: \sum_{\sigma\leq^{\!R}\tau}\sigma\otimes\tau/^{\!R}\sigma
\end{equation}
to describe the map $\delta_{\!R}$. The sum is implicitly indexed by elements $\sigma\in\mcT$ in the canonical basis $\mcB$ of $\mcT$.

We need to introduce  some basic definitions/results on the paraproduct ${\sf P}$ used in the representation Theorem \ref{ThmBH1} before stating our main result. Recall from \cite{BH1} the definition of the two-parameter extension $\overline{\sf P}$ of the paraproduct operator in terms of the kernels $Q_i$ of the Littlewood-Paley projectors -- see e.g. Section 3.1 of \cite{BH1}. For $j\geq 1$, set 
$$
P_j := \sum_{-1\leq i\leq j-2} Q_i,
$$
and for a two variables real-valued distribution $\Lambda$ on $\bbR^d\times\bbR^d$, and $j\geq 1$, set for all $x\in\bbR^d$
$$
\big({\bf Q}_j\Lambda\big)(x) := \big\langle \Lambda , P_j(x-\cdot)\otimes Q_j(x-\cdot)\big\rangle.
$$
The action of $\overline{\sf P}$ on $\Lambda$ is given by 
$$
\overline{\sf P}\Lambda := \sum_{j\geq 1} {\bf Q}_j\Lambda.
$$
It coincides with the paraproduct operator when applied to product distrutions $\Lambda_{y,z} = a(y)b(z)$, in the sense that
$$
\overline{\sf P}\big(a(y)b(z)\big) = {\sf P}_ab.
$$
We use here a formal notation to emphasize the dependence of a distribution on $\bbR^d\times\bbR^d$ on its arguments. Recall also from \cite{BB3} the definition of the operator
\begin{equation} \label{Commutator}
{\sf R}(a,b,c) := {\sf P}_a({\sf P}_bc) - {\sf P}_{ab}c,
\end{equation}
and the fact that it maps continuously $C^\alpha(\bbR^d)\times C^\beta(\bbR^d)\times C^\gamma(\bbR^d)$ into $C^{\alpha+\beta+\gamma}(\bbR^d)$, for all $\alpha,\beta\in (0,1)$ and $\gamma\in(-3,3)$. (See Proposition 3 in \cite{BB3} -- the parameters in the definition of the operators can be arranged so as to get the continuity of $\sf R$ for $\gamma$ in any a priori fixed interval of regularity exponents. The interval $(-3,3)$ has thus no special meaning.) We need also a key recursive identity which has been used in \cite{BH1} -- identity (2.5) therein. Rewriting the identity
\begin{equation*}
{\sf \Pi}^{\!R}\tau =   \sum_{\sigma\leq \tau} {\sf g}^{\!R}_x(\tau/\sigma) \, {\sf \Pi}^{\!R}_x\sigma
\end{equation*}
under the form
\begin{equation*}
{\sf \Pi}_x^{\!R} \tau = {\sf \Pi}^{\!R}\tau -  \sum_{\sigma < \tau} {\sf g}^{\!R}_x(\tau/\sigma) \, {\sf \Pi}^{\!R}_x\sigma 
\end{equation*}
and iterating, we get first
\begin{equation*}
{\sf \Pi}_x^{\!R} \tau = {\sf \Pi}^{\!R}\tau -  \sum_{\sigma < \tau} {\sf g}^{\!R}_x(\tau/\sigma) \, {\sf \Pi}^{\!R} \sigma + \sum_{\sigma_2 < \sigma_1 < \tau} {\sf g}^{\!R}_x(\tau/\sigma_1) {\sf g}^{\!R}_x(\sigma_1/\sigma_2)\,{\sf \Pi}^{\!R}_x \sigma_2, 
\end{equation*}
and after a finite number of iterations
\begin{equation} \label{recurPi}
{\sf \Pi}_x^{\!R} \tau = {\sf \Pi}^{\!R}\tau -  \sum_{n \geq 1}(-1)^n\ \sum_{\sigma_n < \cdots <\sigma_1 < \tau} {\sf g}^{\!R}_x(\tau/\sigma_1) \cdots {\sf g}^{\!R}_x(\sigma_{n-1}/\sigma_{n}) \, {\sf \Pi}^{\!R} \sigma_n.
\end{equation}
Similarly, one has 
\begin{equation} \label{recurPi2}
{\sf \Pi}_x \tau = {\sf \Pi}\tau -  \sum_{n \geq 1}(-1)^n\ \sum_{\sigma_n < \cdots <\sigma_1 < \tau} {\sf g}_x(\tau/\sigma_1) \cdots {\sf g}_x(\sigma_{n-1}/\sigma_{n}) \, {\sf \Pi}\sigma_n.
\end{equation}
If one uses relation \eqref{EqPiMx} to write
$$
{\sf \Pi}^{\!R}_x\tau  = \sum_{\sigma\leq^{\!R}\tau} {\sf g}^{-1}_x(\tau/^{\!R}\sigma) \, {\sf \Pi}_x\sigma
$$
we obtain from \eqref{recurPi2} the identity
\begin{equation} \label{EqIdentityPixM}
{\sf \Pi}^{\!R}_x\tau = \sum_{\sigma\leq^{\!R}\tau} {\sf g}^{-1}_x(\tau/^{\!R}\sigma) \, {\sf \Pi} \sigma - 
\sum_{n \geq 1}(-1)^n\ \sum_{\sigma_n < \cdots <\sigma_1 < \sigma \leq^{\!R} \tau} {\sf g}^{-1}_x(\tau/^{\!R}\sigma) {\sf g}_x(\sigma /\sigma_1) \cdots {\sf g}_x(\sigma_{n-1}/\sigma_{n}) \, {\sf \Pi} \sigma_n.
\end{equation}

\medskip

\begin{thm} \label{ThmMain} 
The formula
\begin{equation} \label{main_result} \begin{aligned}
\, [\tau ]^{\!R}  = & \sum_{n\geq 1}(-1)^{n-1}\sum_{\textsf{\textbf{1}}<\tau_{n+1}<\cdots<\tau_1<\tau} {\sf R}\Big({\sf g}^{\!R}(\tau/\tau_1)\cdots{\sf g}^{\!R}(\tau_{n-1}/\tau_n)\,;\,{\sf g}^{\!R}(\tau_n/\tau_{n+1})\,;\,[\tau_{n+1}]^{\!R}\Big) \\ & + \, \overline{\sf P}\Big(\big({\sf \Pi}_y^{\!R}\tau\big)(z)\Big)+ {\sf S} ( {\sf \Pi}^{\!R}\tau),
\end{aligned} \end{equation}
where $ \overline{\sf P}\Big(\big({\sf \Pi}_y^{\!R}\tau\big)(z)\Big)$ is given by
\begin{equation} \label{EqInductionRelation2} \begin{split}
\overline{\sf P}\Big(\big({\sf \Pi}_y^{\!R}\tau\big)(z)\Big) = \sum_{\sigma\leq^{\!R} \tau} &{\sf P}_{{\sf g}^{-1}(\tau/^{\!R}\sigma)} [\sigma]    + \sum_{n \geq 1} \, \sum_{\textsf{\textbf{1}}<\sigma_n<\cdots<\sigma_1<\sigma \leq^{\!R} \tau} (-1)^{n-1}\times   \\
& {\sf R}\Big({\sf g}^{-1}(\tau/^{\!R}\sigma)\,{\sf g}(\sigma/\sigma_1)\cdots{\sf g}(\sigma_{n-1}/\sigma_n)\,;\,{\sf g}(\sigma_m/\sigma_{n+1})\,;\,[\sigma_{m+1}]\Big),
\end{split}  \end{equation}
and defines inductively the bracket map $[\,\cdot\,]^{\!R}$ in terms of the bracket map $[\,\cdot\,]$. One has moreover
\begin{align*}
{\sf P}_{1} \left({\sf \Pi}^{\!R}\tau \right) =: {\sf \Pi}^{\!R}\tau - {\sf S} ( {\sf \Pi}^{\!R}\tau),
\end{align*}
where $ {\sf S} ( {\sf \Pi}^{\!R}\tau ) $ is a smooth term depending continuously in any Hölder topology on the distribution ${\sf \Pi}^{\!R}\tau$.
\end{thm}

\medskip

\begin{Dem}
One can repeat safely part of the proof of Proposition 12 in \cite{BH1}. We proceed by induction on the size of the trees.
 By applying $ {\sf P} $ to the identity \eqref{recurPi}, one has 
$$
{\sf P}_{1} \left({\sf \Pi}^{\!R}\tau \right)= \sum_{n\geq 1} (-1)^n\sum_{\tau_n<\cdots<\tau_1<\tau} {\sf P}_{{\sf g}^{\!R}(\tau/\tau_1)\dots{\sf g}^{\!R}(\tau_{n-1}/\tau_n)} {\sf \Pi}^{\!R}\tau_n + \overline{\sf P}\Big(\big({\sf \Pi}_y^{\!R}\tau\big)(z)\Big),
$$
 In the end, we have

\begin{align*}
{\sf \Pi}^{\!R}\tau = \sum_{n\geq 1} (-1)^n\sum_{\textsf{\textbf{1}}< \tau_n<\cdots<\tau_1<\tau} {\sf P}_{{\sf g}^{\!R}(\tau/\tau_1)\dots{\sf g}^{\!R}(\tau_{n-1}/\tau_n)} {\sf \Pi}^{\!R}\tau_n + \overline{\sf P}\Big(\big({\sf \Pi}_y^{\!R}\tau\big)(z)\Big) + {\sf S} ( {\sf \Pi}^{\!R}\tau)
\end{align*}
We replace $ \tau_n $ by the following expression
$$
{\sf \Pi}^{\!R}\tau_n = \sum_{\textsf{\textbf{1}}<\tau_{n+1}<\tau_n} {\sf P}_{{\sf g}^{\!R}(\tau_n/\tau_{n+1})}[\tau_{n+1}]^{\!R} + [\tau_n]^{\!R},
$$ 
and using the definition \eqref{Commutator} of the operator $\sf R$, we get 
\begin{equation} \label{EqInductionRelation1} \begin{aligned}
\, & [\tau]^{\!R} = \sum_{n\geq 1}(-1)^{n-1}\sum_{\textsf{\textbf{1}} <\tau_{n+1}<\cdots<\tau_1<\tau} {\sf R}\Big({\sf g}^{\!R}(\tau/\tau_1)\cdots{\sf g}^{\!R}(\tau_{n-1}/\tau_n)\,;\,{\sf g}^{\!R}(\tau_n/\tau_{n+1})\,;\,[\tau_{n+1}]^{\!R}\Big)    \\ 
&+ \, \overline{\sf P}\Big(\big({\sf \Pi}_y^{\!R}\tau\big)(z)\Big) + {\sf S} ( {\sf \Pi}^{\!R}\tau),
\end{aligned} \end{equation}
from the same `fantastic' telescopic sum as in the proof of Proposition 12 in \cite{BH1}. The same mechanics is at work in the proof of identity \eqref{EqInductionRelation2}. Indeed, since one has from identity \eqref{EqIdentityPixM}
\begin{equation*}  \begin{split}
& \overline{\sf P}\Big(\big({\sf \Pi}_y^{\!R}\tau\big)(z)\Big) = \sum_{\sigma\leq^{\!R} \tau} {\sf P}_{{\sf g}^{-1}(\tau/^{\!R}\sigma)} [\sigma] - \\ & \sum_{n\geq 1} (-1)^m\sum_{\sigma_n<\cdots<\sigma_1<\sigma \leq^{\!R} \tau} {\sf P}_{{\sf g}^{-1}(\tau/^{\!R}\sigma)\,{\sf g}(\sigma/\sigma_1)\cdots{\sf g}(\sigma_{n-1}/\sigma_n)}[\sigma_n],
\end{split}  \end{equation*}
and 
$$
{\sf \Pi}\sigma_n = \sum_{\sigma_{n+1}\leq \sigma_n} {\sf P}_{{\sf g}(\sigma_n/\sigma_{n+1})} [\sigma_{n+1}],
$$
a telescopic sum appears and leaves formula \eqref{EqInductionRelation2}. Formulas \eqref{EqInductionRelation1} and \eqref{EqInductionRelation2} give jointly an inductive formula giving $[\tau]^{\!R}$ in terms of the $[\tau']$, with $\tau'\in \mcT$. 
\end{Dem}

\medskip

\begin{cor} \label{CorModelT}
The model $\big(g^{\!R},{\sf \Pi}^{\!R}\big)$ on $\mathscr{T}^{(\epsilon)}$ is actually a model on $\mathscr{T}$.
\end{cor}

\medskip

\begin{Dem}
Given $\tau\in T$ with $\textsf{deg}(\tau)\leq 0$, we know from Proposition 10 in \cite{BH1} that the double sum in \eqref{main_result} defines an element of $C^{\textsf{deg}(\tau)}(\bbR^d)$. Since the distribution $\Lambda=\big({\sf \Pi}_y^{\!R}\tau\big)(z)$ on $\bbR^d_y\times\bbR^d_z$ satisfies from \eqref{EqPiMx} the estimate
$$
\big\|{\bf Q}_j\Lambda\big\| \lesssim 2^{-j\textsf{deg}(\tau)},
$$
uniformly in $j\geq 1$, Proposition 8 in \cite{BH1} tells us that $\overline{\sf P}\Big(\big({\sf \Pi}_y^{\!R}\tau\big)(z)\Big)$ is also an element of $C^{\textsf{deg}(\tau)}(\bbR^d)$. All the brackets $[\tau]^{\!R}$ are thus elements of $C^{\textsf{deg}(\tau)}(\bbR^d)$, so $\big(g^{\!R},{\sf \Pi}^{\!R}\big)$ turns out to be a model on $\mathscr{T}$ from Theorem \ref{ThmBH1}, as the unique model on $\mathscr{T}$ associated to the brackets $[\,\cdot\,]^{\!R}$ provides canonically a model on $\mathscr{T}^{(\epsilon)}$ that needs to coincide with $\big(g^{\!R},{\sf \Pi}^{\!R}\big)$, by uniqueness.
\end{Dem}

\bigskip
\bigskip

\noindent \textcolor{gray}{$\bullet$} {\sf I. Bailleul} -- Univ. Rennes, CNRS, IRMAR - UMR 6625, F-35000 Rennes, France.   \\
\noindent {\it E-mail}: ismael.bailleul@univ-rennes1.fr   

\medskip

\noindent \textcolor{gray}{$\bullet$} {\sf Y. Bruned} --  School of Mathematics, University of Edinburgh, EH9 3FD, Scotland.   \\
{\it E-mail}: Yvain.Bruned@ed.ac.uk


\bigskip
\bigskip



\begin{thebibliography}{99}

\bibitem{BailleulCassWeidner}
I. Bailleul,
\newblock {\em On the definition of a solution to a rough differential equation}.
\newblock Ann. Fac. Sci. Toulouse, {\bf 30}(3):463--478, (2021).

\bibitem{BB3}
I. Bailleul and F. Bernicot,
\newblock {\em High order paracontrolled calculus}.
\newblock Forum of Mathematics Sigma, {\bf 7}(e44):1--93, (2019).

\bibitem{BailleulBruned}
I. Bailleul and Y. Bruned,
\newblock {\em Renormalised singular stochastic PDEs}. 
\newblock arXiv:2101.11949, (2021).

\bibitem{BH1}
I. Bailleul and M. Hoshino,
\newblock {\em Paracontrolled calculus and regularity structures I}.
\newblock J. Math. Soc. Japan, DOI: 10.2969/jmsj/81878187:1--43, (2020).

\bibitem{BH2}
I. Bailleul and M. Hoshino,
\newblock {\em Paracontrolled calculus and regularity structures II}.
\newblock To appear in J. \'Ec. Polytechnique, 1--40, (2021).

\bibitem{RSGuide}
I. Bailleul and M. Hoshino,
\newblock {\em A tourist guide to regularity structures and singular stochastic PDEs}.
\newblock arXiv:2006:03524, 1--81, (2020).

\bibitem{BFPP}
C. Bellingeri and P. Friz and S. Paycha and R. Preiss,
\newblock {\em Smooth rough paths, their geometry and algebraic renormalization}
\newblock arXiv:2111.15539, (2021).

\bibitem{bruned:tel-01306427}
Y.~Bruned.
\newblock \emph{{Singular KPZ Type Equations}}.
\newblock 205 pages, {PhD thesis}, {Universit{\'e} Pierre et Marie Curie - Paris VI},
  2015.
\newblock \url{https://tel.archives-ouvertes.fr/tel-01306427}.

\bibitem{BrunedRecursive}
Y. Bruned,
\newblock {\em Recursive formulae for regularity structures}.
\newblock Stoch. PDEs: Anal. Comp., {\bf 6}(4):525--564, (2018).

\bibitem{BrunedTZ}
Y. Bruned,
\newblock {\em Renormalization from non-geometric to geometric rough paths}.
\newblock To appear in Ann. Inst. H. Poincar\'e Probab. Statist.,  arXiv:2007.14385,1--19, (2020).

\bibitem{BCCH18}
Y. Bruned and A. Chandra and I. Chevyrev and M. Hairer,
\newblock {\em Renormalising SPDEs in regularity structures}.
\newblock J. Europ. Math. Soc., {\bf 23}(3):869--947, (2021).

\bibitem{BHZ}
Y. Bruned and M. Hairer, and L. Zambotti,
\newblock {\em Algebraic renormalization of regularity structures},
\newblock Invent. Math., {\bf 215}(3):1039--1156, (2019).

\bibitem{BrunedManchon}
Y. Bruned and D. Manchon,
\newblock {\em Algebraic deformation for (S)PDEs}. 
\newblock arXiv:2011.05907, (2020).

\bibitem{CassWeidner}
T. Cass and M. Weidner,
\newblock {\em Tree algebras over topological vector spaces in rough path theory}.
\newblock  arXiv:1604.07352, (2016).

\bibitem{ChandraHairer}
A. Chandra and M. Hairer,
\newblock {\em An analytic BPHZ theorem for Regularity Structures}.
\newblock arXiv:1612.08138, (2016).

\bibitem{ChandraWeber}
A. Chandra and H. Weber,
\newblock {\em Stochastic PDEs, regularity structures and interacting particle systems}.
\newblock Ann. Fac. Sci. Toulouse, {\bf 26}(4):847--909, (2017).

\bibitem{CorwinShen}
I. Corwin and H. Shen,
\newblock {\em Some recent progress in singular stochastic PDEs}.
\newblock Bull. Am. Math. Soc., {\bf 57}:409--454, (2020).

\bibitem{FrizHairer}
P. Friz and M. Hairer,
\newblock {\em A course on rough paths, with an introduction to regularity structures}
\newblock Universitext, Springer, (2020).

\bibitem{GubinelliBranched}
M. Gubinelli,
\newblock {\em Ramification of rough paths}.
\newblock J. Diff. Eq., {\bf 248}(4):693--721, (2010).

\bibitem{Hai14}
M. Hairer,
\newblock {\em A theory of regularity structures}.
\newblock Invent. Math., {\bf 198}(2):269--504, (2014).

\bibitem{HairerBrazil}
M. Hairer,
\newblock {\em Introduction to Regularity Structures}.
\newblock Braz. Jour. Prob. Stat., {\bf 29}(2):175--210, (2015).

\bibitem{HairerTakagi}
M. Hairer,
\newblock {\em renormalization of parabolic stochastic PDEs}.
\newblock Japanese. J. Math., {\bf 13}:187--233, (2018).

\bibitem{HairerKelly}
M. Hairer and D. Kelly,
\newblock {\em Geometric versus non-geometric rough paths}.
\newblock Ann. Institut H. Poincar\'e, {\bf 51}(1):207--251, (2015).

\bibitem{HQ}
M. Hairer and J. Quastel,
\newblock {\em A class of growth models rescaling to KPZ}.
\newblock Forum Math. Pi, {\bf 6}(e3):1--112 , (2018).

\bibitem{SinghTeichmann}
H. Singh and J. Teichmann,
\newblock {\em An elementary proof of the reconstruction theorem},
\newblock arXiv:1802.03082, (2018).

\bibitem{TZ18}
N. Tapia and L. Zambotti,
\newblock {\em The geometry of the space of branched Rough Paths},
\newblock Proc. London Math. Soc., {\bf 121}(2):220-251, (2020).

\end{thebibliography}
\end{document}